\documentstyle[12pt,amssymb,amsbsy,righttag]{amsart}

\renewcommand{\a}{\alpha}
\renewcommand{\b}{\beta}
\newcommand{\g}{\gamma}
\renewcommand{\d}{\delta}
\newcommand{\e}{\varepsilon}
\newcommand{\z}{\zeta}

\renewcommand{\l}{\lambda}
\newcommand{\r}{\rho}
\newcommand{\s}{\sigma}
\renewcommand{\t}{\tau}

\newcommand{\f}{\varphi}
\renewcommand{\o}{\omega}

\newcommand{\G}{\Gamma}

\renewcommand{\L}{\Lambda}

\renewcommand{\O}{\Omega}

\newcommand{\B}{{\cal B}}

\newcommand{\M}{{\cal M}}

\newcommand{\1}{{\bf 1}}

\newcommand{\C}{{\Bbb C}}
\newcommand{\T}{{\Bbb T}}
\newcommand{\pp}{{\Bbb P}}
\newcommand{\dd}{{\Bbb D}}
\newcommand{\R}{{\Bbb R}}
\newcommand{\Z}{{\Bbb Z}}

\newcommand{\0}{{\Bbb O}}

\newcommand{\bs}{\boldsymbol}

\newcommand{\m}{{\boldsymbol m}}
\newcommand{\bS}{{\boldsymbol S}}

\newcommand{\rf}[1]{(\ref{#1})}

\newcommand{\df}{\stackrel{\mathrm{def}}{=}}
\newcommand{\dist}{\operatorname{dist}}

\newcommand{\re}{\operatorname{Re}}

\newcommand{\supp}{\operatorname{supp}}

\newcommand{\rank}{\operatorname{rank}}
\newcommand{\const}{\operatorname{const}}

\newcommand{\eeq}{\end{equation}}
\newcommand{\beq}{\begin{equation}}
\newcommand{\bay}{\begin{eqnarray}}
\newcommand{\ey}{\end{eqnarray}}
\newcommand{\bey}{\begin{eqnarray*}}
\newcommand{\eey}{\end{eqnarray*}}

\newcommand{\be}{\infty}

\newcommand{\bl}{\blacksquare}

\newcommand{\Pf}{{\bf Proof. }}

\renewcommand{\re}{\operatorname{Re}}
\newcommand{\ov}{\overline}

\newtheorem{thm}{\hspace{\parindent}Theorem}[section]

\newtheorem{cor}[thm]{\hspace{\parindent}Corollary}
\newtheorem{lem}[thm]{\hspace{\parindent}Lemma}

\begin{document}

\newcommand{\Mp}{{\frak M}_p}
\newcommand{\card}{\operatorname{card}}
\newcommand{\Mo}{{\frak M}_p^0}
\newcommand{\st}{\star}
\newcommand{\fl}{\flat}
\newcommand{\PR}{\hat{\otimes}}
\renewcommand{\theequation}{\thesection.\arabic{equation}}

\author{A.B. Aleksandrov and V.V. Peller}
\thanks{The first author is partially supported by Grant 99-01-00103
of Russian Foundation of Fundamental Studies and by Grant 326.53 of Integration.
The second author is partially supported by NSF grant DMS 9970561.}

\title{Hankel and Toeplitz--Schur Multipliers}

\maketitle

\begin{abstract}
We study the problem of characterizing Hankel--Schur multipliers and 
Toeplitz--Schur multipliers of Schatten--von Neumann class $\bS_p$ for 
\linebreak$0<p<1$. We obtain various sharp necessary conditions and 
sufficient conditions for a Hankel matrix 
to be a Schur multiplier of $\bS_p$. We also give a
characterization of the Hankel--Schur multipliers of $\bS_p$ whose symbols
have lacunary power series.
Then the results on Hankel--Schur multipliers are used to obtain
a characterization of the Toeplitz--Schur multipliers of $\bS_p$.
Finally, we return to Hankel--Schur multipliers and obtain new results
in the case when the symbol of the Hankel matrix is a complex measure on
the unit circle.
\end{abstract}

\setcounter{equation}{0}
\setcounter{section}{0}
\section{\bf Introduction}

\

The {\it Schur product} of matrices $A=\{a_{jk}\}_{j,k\ge0}$ and 
$B=\{b_{jk}\}_{j,k\ge0}$ is defined as the matrix $A\star B$ whose entries
are the products of the entries of $A$ and $B$:
$$
A\star B=\{a_{jk}b_{jk}\}_{j,k\ge0}.
$$
If we identify in a natural way the bounded linear operators on $\ell^2$ 
with their matrix representation with respect to the standard orthonormal basis of
$\ell^2$, we can study Schur multipliers of various classes of linear 
operators on $\ell^2$. Namely, if $\cal X$ is a class of bounded linear operators
on $\ell^2$, we say that a matrix $A$ is a {\it Schur multiplier} of $\cal X$ if
and only if 
$$
B\in \cal X\Longrightarrow A\star B\in \cal X.
$$
We are interested in this paper in Schur multipliers of Schatten--von Neumann
classes $\bS_p$ (see [GK], [BS1] for information on the classes $\bS_p$). For $0<p<\be$
we denote by ${\frak M}_p$ the space of Schur multipliers of $\bS_p$ and we put
$$
\|A\|_{\Mp}=\sup\{\|A\st B\|_{\bS_p}:~\|B\|_{\bS_p}\le1\}.
$$
It is easy to see that for $p\ge1$ the functional $\|\cdot\|_{\Mp}$ is a norm
on $\Mp$. Using the triangle inequality
\beq
\label{1.0}
\|T_1+T_2\|_{\bS_p}^p\le\|T_1\|_{\bS_p}^p+\|T_2\|_{\bS_p}^p,\quad 0<p\le1,
\end{equation}
(see [Pie], [BS1]), one can easily see 
that 
\beq
\label{1.1}
\|A_1+A_2\|^p_{\Mp}\le\|A_1\|^p_{\Mp}+\|A_2\|^p_{\Mp},\quad0<p\le1.
\end{equation}
We denote by $\frak M$ the class of Schur multipliers of the space $\B$ of bounded
linear operators on $\ell^2$.

In this paper we are going to study the {\it Hankel--Schur multipliers} of
$\bS_p$, i.e., matrices of class $\Mp$ of the form $\{\g_{j+k}\}_{j,k\ge0}$ 
and the {\it Toeplitz--Schur multipliers} of
$\bS_p$, i.e., matrices of class $\Mp$ of the form $\{t_{j-k}\}_{j,k\ge0}$.

Let us summarize briefly some well-known properties of classes $\Mp$. The class
${\frak M}_2$ is the space of matrices with bounded entries. If $1<p<\be$, then
$\Mp={\frak M}_{p'}$, where $1/p+1/p'=1$. The space $\frak M$ coincides with 
${\frak M}_1$. This follows from the facts that the dual space to $\bS_p$ can
be identified in a natural way with $\bS_{p'}$ and the dual to $\bS_1$ can
be identified with $\B$.

Next, interpolating between the classes $\bS_p$, one can easily see that 
$$
{\frak M}_{p_1}\subset{\frak M}_{p_2}\quad\mbox{if}\quad0<p_1\le p_2\le2.
$$

To describe the space $\frak M$, we consider the projective tensor product 
$\ell^\be\PR\ell^\be$ that consists
of matrices \mbox{$C=\{c_{jk}\}_{j,k\ge0}$} for which there exist sequences
\linebreak$X^{(n)}=\left\{x^{(n)}_j\right\}_{j\ge0}\in\ell^\be$ and 
$Y^{(n)}=\left\{y^{(n)}_j\right\}_{j\ge0}\in\ell^\be$, $n\ge0$, such that
\beq
\label{1.2}
c_{jk}=\sum_{n\ge0}x^{(n)}_jy^{(n)}_k
\end{equation}
and
\beq
\label{1.3}
\sum_{n\ge0}\left\|X^{(n)}\right\|_{\ell^\be}\left\|Y^{(n)}\right\|_{\ell^\be}<\be.
\end{equation}
The norm $\|C\|_{\ell^\be\PR\ell^\be}$ is by definition the {\it infimum} in
\rf{1.3} over all sequences $X^{(n)}$ and $Y^{(n)}$ satisfying \rf{1.2}.

For any positive integer $m$ we consider the matrix 
$Q^{(m)}=\{q^{(m)}_{jk}\}_{j,k\ge0}$ defined by
\beq
\label{1.4}
q^{(m)}_{jk}=
\left\{\begin{array}{ll}1,&j\le m,~k\le m,\\0,&\mbox{otherwise}\,.\end{array}\right.
\end{equation}
It is well known that a matrix $A$ belongs to ${\frak M}={\frak M}_1$ if and only
if
$$
\sup_m\|Q^{(m)}\st A\|_{\ell^\be\PR\ell^\be}<\be
$$
(see [Be]).

One can generalize the notion of a Schur multiplier of $\bS_p$ to the
case of linear operators on $L^2(X,\mu)$, where $\mu$ 
is a Borel measure on a metric space $X$. For $p\le2$ the space
$\Mp(\mu)$ of Schur multipliers of $\bS_p$ is defined as the class of 
measurable functions $\f$ on $X\times X$ such that for any integral operator 
$T\in\bS_p$ with kernel function $k$, the integral operator with kernel
function $\f k$ is also in $\bS_p$. For \linebreak$2<p<\be$ we put 
$\Mp(\mu)\df{\frak M}_{p'}(\mu)$ and the space ${\frak M}(\mu)$ of Schur multipliers
of the space of bounded linear operators on $L^2(\mu)$ is by definition
${\frak M}_1(\mu)$.

Note that Schur multipliers are a powerful tool to investigate linear operators,
they are applied in perturbation theory, in operator algebras, operator spaces,
in the study of Hankel operators, and many other fields of mathematics,
see e.g., [Be], [BS2], [Pel2], [Pel3], [Pel4], [Pis1], [Pis2].

The Hankel--Schur multipliers and the Toeplitz--Schur multipliers of $\bS_p$ are well 
investigated for $p=1$ (see \S 2 for a brief summary of the results). In this
paper we consider the case $p<1$.

In \S 4 we study Hankel--Schur multipliers of $\bS_p$, $0<p<1$. We find many
sharp necessary conditions and sufficient conditions. We also obtain in \S 4 a necessary and 
sufficient condition for a Hankel matrix with a ``lacunary'' symbol to be a 
Schur multiplier of $\bS_p$. In \S 5 we use the results of \S 4 to 
study the Toeplitz--Schur multipliers of $\bS_p$ and we obtain their characterization 
which can be formulated very easily. Finally, in \S 6 we return to Hankel--Schur multipliers.
Namely, we will study the case when the symbol of a Hankel matrix is a complex measure on the 
unit circle. Note that the results obtained in \S 6 are based on the results of
\S 5 on Toeplitz--Schur multipliers. In \S 2 we collect necessary information 
and in \S 3 we obtain some useful general results on Schur multipliers of $\bS_p$.

\

\setcounter{equation}{0}
\setcounter{section}{1}
\section{\bf Preliminaries}

\

{\bf $\bs{L^p}$-norms of certain trigonometric polynomials.} Let $F$ be an 
infinitely differentiable function on $\R$
with compact support and let $m$ be a positive integer. Consider the trigonometric
polynomial
\bay
\label{2.z}
F_{(m)}(z)=\sum_{k\in\Z}F\left(\frac{k}{m}\right)z^k.
\ey
We need the following fact:
\beq
\label{2.1}
c_F\, m^{1-1/p}\le\|F_{(m)}\|_{L^p(\T)}\le C_F\, m^{1-1/p},\quad0<p<\be.
\end{equation} 
for some constants $c_F$ and $C_F$. We refer the reader to [A1] for the proof.

{\bf Hankel operators of class $\bs{\bS_p}$.} Let $\f=\sum_{j\ge0}\hat\f(j)z^j$ 
be a function analytic in the unit disk $\dd$.
We can associate with $\f$ the {\it Hankel matrix} 
$\G_\f\df\{\hat\f_{j+k}\}_{j,k\ge0}$. If $\f$ belongs to the Hardy class $H^2$,
we can consider the matrix $\G_\f$ as an operator in $\ell^2$ whose domain
is the set of finitely supported sequences. We call such operators
{\it Hankel operators}. By Nehari's theorem (see e.g., [N]),
$\G_\f$ is a bounded operator on $\ell^2$ if and only if there exists a
function $\psi$ in $L^\be$ whose Fourier coefficients satisfy
$$
\hat\psi(j)=\hat\f(j),\quad j\in\Z_+.
$$
The following result describes the Hankel operators of class $\bS_p$:
\beq
\label{2.2}
\G_\f\in\bS_p\quad\mbox{\it if and only if}\quad\f\in B_p^{1/p},\quad0<p<\be,
\end{equation}
where $B_p^{1/p}$ is a Besov class which will be defined below. For $p>1$
\rf{2.2} was proved in [Pel1] (the implication $\Leftarrow$ in the
case $0<p<1$ was also established in [Pel1]), while for $0<p<1$ the result
was been obtained independently in [Pel3], [Pek], and [Se].

We can also define a Hankel matrix $\G_\f$ for any distribution $\f$ on $\T$:
\linebreak$\G_\f=\{\hat\f_{j+k}\}_{j,k\ge0}$. Clearly, $\G_\f=\G_{\pp_+\f}$, where
$\pp_+\f\df\sum_{n\ge0}\hat\f(n)z^n$.

The {\it Besov spaces} $B_{pq}^s$, $0<p\le\be$, $0<q\le\be$, $s\in\R$, 
of functions (or distributions) on $\T$ can be defined in
many different ways (see [Pee]). We give a definition most suitable for our 
purposes. Let $v$ be a function in $C^\be(\R)$ such that $v\ge0$, 
$\supp v=[1/2,2]$, and $\sum_{n\ge0}v(x/2^n)=1$ for $x\ge1$. The 
trigonometric polynomials $V_n$ are defined by
\bay
\label{bes}
V_n=\sum_{k\ge0}v\left(\frac{k}{2^n}\right)z^k, \quad n\ge1,
\quad V_0(z)=\bar z+1+z,\quad V_n=\ov{V_{-n}},\quad n<0.
\ey
The Besov class $B_{pq}^s$, $0<p\le\be$, $0<q\le\be$, $s\in\R$, consists of 
distributions $\psi$ on $\T$ such that
$$
\left\{2^{|n|s}\|\psi*V_n\|_p\right\}_{n\in\Z}\df
\left\{2^{|n|s}\|\psi*V_n\|_{L^p}\right\}_{n\in\Z}
\in\ell^q(\Z).
$$
The classes $B_{pq}^s$ do not depend on the choice of a particular function $v$.

We denote by $\left(B_{pq}^s\right)_+$ the subspace of $B_{pq}^s$ that consists of
functions analytic in $\dd$, i.e., a function $\psi$ analytic in $\dd$ belongs to 
$\left(B_{pq}^s\right)_+$ if and only if
$$
\left\{2^{ns}\|\psi*V_n\|_{L^p}\right\}_{n\ge0}\in\ell^q.
$$
We use the notation $B_p^s$ for $B_{pp}^s$.

In particular, the most important Besov class $B_p^{1/p}$, 
$0<p<\be$, consists of distributions $\psi$ on $\T$ such that
$$
\sum_{n\in\Z}2^{|n|}\|\psi*V_n\|_p^p<\be.
$$

{\bf Hankel and Toeplitz--Schur multipliers of $\bs{\bS_1}$.} The description
of Toeplitz--Schur multipliers is well-known, though we do not know who did it first.
Perhaps this is a mathematical folklore. Anyway, a Toeplitz matrix 
$\{t_{j-k}\}_{j,k\ge0}$ is a Schur multiplier of $\bS_1$ if and only if
there exists a complex Borel measure $\mu$ on $\T$ such that $t_j=\hat\mu(j)$,
$j\in\Z$, (we use the notation $\hat\mu(j)$ for the $j$th Fourier coefficients of $\mu$).
This fact can be proved very easily. First, if $\mu=\d_\t$, the point mass
at $\t\in\T$, then 
\beq
\label{2.3}
\hat\mu(j-k)=\bar\t^j\t^k,
\end{equation}
and so $\|\{\mu(j-k)\}_{j,k\ge0}\|_{\ell^\be\PR\ell^\be}\le1$ which proves
that $\{\mu(j-k)\}_{j,k\ge0}\in{\frak M}_1$. Moreover, it is easy to see from
\rf{2.3} that
$\|\{\mu(j-k)\}_{j,k\ge0}\|_{{\frak M}_1}=1$. The fact that 
$\{\mu(j-k)\}_{j,k\ge0}\in{\frak M}_1$ for an arbitrary $\mu$ follows immediately from
the fact that any $\mu$ in the unit ball of the space of Borel measures can be
approximated in the weak-* topology by measures of the form
$$
\sum_j c_j\d_{\t_j},\quad \t_j\in\T,\quad\sum_j|c_j|\le1.
$$

Conversely, if $\{t_{j-k}\}_{j,k\ge0}\in{\frak M}_1=\frak M$, then for any bounded
Toeplitz matrix $\{\a_{j-k}\}_{j,k\ge0}$ the Toeplitz matrix 
$\{t_{j-k}\a_{j-k}\}_{j,k\ge0}$ is also the matrix of a bounded Toeplitz operator.
Since $\{\a_{j-k}\}_{j,k\ge0}\in\B$ if and only if $\a_j=\hat g(j)$ for a function
$g$ in $L^\be$ (see e.g., [D]), it follows that $\{t_{j-k}\}_{j,k\ge0}\in\frak M$
if and only if the sequence $\{t_j\}_{j\in\Z}$ is a Fourier multiplier of the
space $L^\be$, i.e.,
\beq
\label{2.4}
\sum_{j\in\Z}\hat g(j)z^j\in L^\be\quad\Longrightarrow\quad
\sum_{j\in\Z}t_j\hat g(j)z^j\in L^\be.
\end{equation}
It is well known (see e.g., [K]) that \rf{2.4} is equivalent to the fact that
$\{t_j\}_{j\in\Z}$  is a sequence of the Fourier coefficients of a complex Borel
measure.

The problem of describing Hankel--Schur multipliers of $\bS_1$ is more complicated.
It was observed in [Pel2] that if a Hankel matrix $\{\g_{j+k}\}_{j,k\ge0}$
is a Schur multiplier of $\bS_1$, then the sequence $\{\g_j\}_{j\ge0}$
is a Fourier multiplier of the Hardy class $H^1$, i.e.,
$$
\sum_{j\in\Z_+}\hat g(j)z^j\in H^1\quad\Longrightarrow\quad
\sum_{j\in\Z_+}\g_j\hat g(j)z^j\in H^1.
$$
A question was asked in [Pel2] of whether the converse is true. It was shown in
[L] that this is not the case. Later a necessary and sufficient condition has been
found in [Pis1]: $\{\g_{j+k}\}_{j,k\ge0}$ is a Schur multiplier of $\bS_1$
if and only if $\{\g_j\}_{j\ge0}$ is a Fourier multiplier of the Hardy space
$H^1(\bS_1)$ of $\bS_1$-valued functions, i.e.,
$$
\sum_{j\in\Z_+}\hat G(j)z^j\in H^1(\bS_1)\quad\Longrightarrow\quad
\sum_{j\in\Z_+}\g_j\hat G(j)z^j\in H^1(\bS_1).
$$
\

\setcounter{equation}{0}
\setcounter{section}{2}
\section{\bf Some General Results on Schur Multipliers}

\

In this section we obtain several general results on Schur multipliers of class 
$\bS_p$ for $0<p<1$. We start with an analogue of a theorem of Schur 
for multipliers of class $\bS_p$ for $0<p<1$.
Recall that by the Schur theorem [Sc] any matrix of a bounded 
linear operator on $\ell^2$ is a Schur multiplier of $\B$. 
The following theorem is an analog of the 
Schur theorem for Schur multipliers of $\bS_p$ in the case $0<p<1$.

Throughout the rest of this paper we associate with a given number $p\in(0,1)$
the positive number ${p_\#}$ defined by
\beq
\label{n.1}
{p_\#}=\frac{p}{1-p}
\end{equation}
It is convenient to assume that if $p=1$, then ${p_\#}=\be$.

We also associate with $p$ the number
\begin{equation}
\label{flat}
p_\flat\df\frac{2p}{2-p}. 
\end{equation}
Clearly, $p_\fl$ makes sense for any $p<2$.

\begin{thm}
\label{tn.1}
Let $0<p<1$ and let ${p_\#}$ be defined by {\em\rf{n.1}}. 
Suppose that $A$ is a matrix of a 
linear operator of class $\bS_{p_\#}$. Then $A\in\Mp$ and
$$
\|A\|_{\Mp}\le\|A\|_{\bS_{p_\#}}.
$$
\end{thm}

\Pf By \rf{1.0}, it suffices to estimate $\|A\st B\|_{\bS_p}$ for rank one operators 
$B$. Let $B=(\cdot,a)b$, where $a$ and $b$ are sequences in $\ell^2$. Consider the 
multiplication operators $D_{\bar a}$ and $D_b$ defined by
$$
D_{\bar a}\{x_n\}_{n\ge0}=\{\bar a_nx_n\}_{n\ge0},\quad
D_b\{x_n\}_{n\ge0}=\{b_nx_n\}_{n\ge0},\quad\{x_n\}_{n\ge0}\in\ell^2.
$$
It is easy to see that
$$
A\st B=D_bAD_{\bar a},
$$
and so
$$
\|A\st B\|_{\bS_p}\le\|D_b\|_{\bS_2}\|A\|_{\bS_{p_\#}}\|D_{\bar a}\|_{\bS_2}
=\|a\|_{\ell^2}\|A\|_{\bS_{p_\#}}\|b\|_{\ell^2}=\|A\|_{\bS_{p_\#}}\|B\|_{\bS_p}
$$
which proves the result. $\bl$

Note that if we take $p=1$, the above reasoning proves that if $A$ is a matrix
of a bounded linear operator on $\ell^2$, then $A$ is a Schur multiplier of $\bS_1$
which is equivalent to the Schur theorem mentioned above.

We are going to study now properties of general block-diagonal and block-triangular
matrices which are Schur multipliers of $\bS_p$ for $0<p\le1$. 
Let $\{m_k\}_{k\ge0}$ and $\{n_l\}_{l\ge0}$ be two increasing sequences of
positive integers such that $m_0=n_0=0$. With each matrix $\O=\{\o_{st}\}_{s,t\ge0}$
we can associate the block-matrix $\{\O_{kl}\}_{k,l\ge0}$ defined by
\beq
\label{n.2}
\O_{kl}=\{\o_{st}\}_{n_k\le s<n_{k+1},\,m_l\le t<m_{l+1}}.
\end{equation}

\begin{thm}
\label{tn.2}
Let $0<p<1$ and let ${p_\#}$ be defined by {\em\rf{n.1}}.
Suppose that the block-matrix $\{\O_{kl}\}_{k,l\ge0}$ defined by {\em\rf{n.2}}
is block-diagonal, i.e., $\O_{kl}=0$ for $k\ne l$. Then
$$
\|\O\|_{\Mp}=\left(\sum_{k\ge0}\|\O_{kk}\|_{\Mp}^{p_\#}\right)^{1/{p_\#}}.
$$
\end{thm}

Note that the corresponding fact for $p=1$ is also valid and well known.

\Pf Let us prove first that 
$\|\O\|^{p_\#}_{\Mp}\le\sum_{k\ge0}\|\O_{kk}\|_{\Mp}^{p_\#}$. 
Let $\{a_s\}_{s\ge0}$ and 
$\{b_t\}_{t\ge0}$ be sequences in $\ell^2$ such that 
$\|a\|_{\ell^2}=\|b\|_{\ell^2}=1$. Then the matrix
$$
\{a_sb_t\}_{s,t\ge0}
$$
has rank one and its $\bS_p$ norm is equal to 1. Clearly, we have to show that
$$
\|\{\o_{st}a_sb_t\}_{s,t\ge0}\|_{\bS_p}^{p_\#}\le\sum_{k\ge0}\|\O_{kk}\|_{\Mp}^{p_\#}.
$$
We have
\bey
\|\{\o_{st}a_sb_t\}_{s,t\ge0}\|_{\bS_p}^p&=&
\sum_{k\ge0}
\left\|\{\o_{st}a_sb_t\}_{m_k\le s<m_{k+1},\,n_k\le t<n_{k+1}}\right\|^p_{\bS_p}\\
&\le&\sum_{k\ge0}\|\O_{kk}\|_{\Mp}^p\left(\sum_{s=m_k}^{m_{k+1}-1}|a_s|^2\right)^{p/2}
\left(\sum_{t=n_k}^{n_{k+1}-1}|b_t|^2\right)^{p/2}\\[.2cm]
&\le&\left(\sum_{k\ge0}\|\O_{kk}\|^{p_\#}_{\Mp}\right)^{1-p}
\eey
by H\"{o}lder's inequality with exponents $1/(1-p)$, $2/p$, and $2/p$.

Let us prove now that
\beq
\label{n.3}
\sum_{k\ge0}\|\O_{kk}\|_{\Mp}^{p_\#}\le\|\O\|^{p_\#}_{\Mp}.
\end{equation}
We can find sequences $\{a_s\}_{s\ge0}$ and $\{b_t\}_{t\ge0}$ such that
$$
\sum_{s=m_k}^{m_{k+1}-1}|a_s|^2=\sum_{t=n_k}^{n_{k+1}-1}|b_t|^2=1
$$
for any $k\ge0$ and
$$
\|\{\o_{st}a_sb_t\}_{m_k\le s<m_{k+1},\,n_k\le t<n_{k+1}}\|_{\bS_p}=\|\O_{kk}\|_{\Mp}.
$$
Let $\{\a_k\}_{k\ge0}$ be a sequence of nonnegative numbers whose norm in $\ell^2$ is
one. Put
$$
\O_k^{(\a)}\df\{\o_{st}\a_k^2a_sb_t\}_{m_k\le s<m_{k+1},\,n_k\le t<n_{k+1}}
$$
and consider the block-diagonal matrix $\O^{(\a)}$ with diagonal entries
$\O_k^{(\a)}$, $k\ge0$. Clearly,
\bey
\|\O\|^p_{\Mp}&\ge&\|\O^{(\a)}\|^p_{\bS_p}
=\sum_{k\ge0}\a_k^{2p}\|\{\o_{st}a_sb_t\}_{s,t\ge0}\|^p_{\bS_p}
=\sum_{k\ge0}\a_k^{2p}\|\O_{kk}\|^p_{\Mp}.
\eey
This inequality holds for an arbitrary sequence $\{\a_k\}_{k\ge0}$ of nonnegative 
numbers whose $\ell^2$ norm is one. Hence,
$$
\sum_{k\ge0}\|\O_{kk}\|^{p_\#}_{\Mp}\le\|\O\|^{p_\#}_{\Mp}.\quad\bl
$$

Interestingly, inequality \rf{n.3} holds not only for block-diagonal matrices but
also for block-triangular ones. To see this, we need the following fact.

\begin{thm}
\label{tn.3}
Suppose that $\O=\{\o_{st}\}_{s,t\ge0}$ is a matrix which has a block-triangular
structure $\{\O_{kl}\}_{k,l\ge0}$, i.e., $\O_{kl}=0$ for $l>k$. Then
$$
\sum\|\O_{kk}\|^p_{\bS_p}\le\|\O\|^p_{\bS_p}.
$$
\end{thm}

\Pf By adding if necessary zero rows or zero columns, we may assume without loss of
generality that $m_k=n_k$ for $k\ge0$. Consider first the case \linebreak
$m_k=n_k=k$, $k\ge0$.
This means that the matrix $\{\o_{st}\}_{s,t\ge0}$ is triangular, i.e., $\o_{st}=0$
for $t>s$. We have to show that 
$$
\sum_{s\ge0}|\o_{ss}|^p\le\|\O\|^p_{\bS_p}.
$$
Put $\O^{(n)}\df\{\o_{st}\}_{0\le s\le n,\,0\le t\le n}$. Note that 
$\{\o_{ss}\}_{0\le s\le n}$ is the sequence of eigenvalues of the matrix $\O^{(n)}$.
Then (see [GK], Ch. 2, \S 3.3)
$$
\sum_{s=0}^n|\o_{ss}|^p\le\|\O^{(n)}\|^p_{\bS_p}.
$$
It remains to make $n$ tend to $\be$.

Suppose now that $m_k$ is an arbitrary increasing sequence of integers in $\Z_+$
such that $m_0=0$ and $n_k=m_k$ for $k\in\Z_+$. Let $U_k$ and $V_k$ be unitary matrices 
of size $(m_{k+1}-m_k)\times(m_{k+1}-m_k)$ such that the matrix $U_k\O_{kk}V_k$ 
is
diagonal. Let $U$ and $V$ be the block diagonal matrices with diagonal entries equal 
to $U_k$, $k\ge0$, and $V_k$, $k\ge0$, respectively. Since we have already proved the
 result for scalar triangular 
matrices, we obtain
$$
\sum_{k\ge0}\|\O_{kk}\|_{\bS_p}^p=\sum_{k\ge0}\|U_k\O_{kk}V_k\|^p_{\bS_p}
\le\|U\O V\|^p_{\bS_p}=\|\O\|_{\bS_p}. \quad\bl
$$

\begin{cor}
\label{tn.4}
Under the hypotheses of Theorem \ref{tn.3} we have
$$
\sum_{k-l=n}\|\O_{kl}\|^p_{\bS_p}\le C(n,p)\|\O\|^p_{\bS_p},
$$
where $C(n,p)$ is a constant that can depend only on $n$ and $p$.
\end{cor}

\Pf The result follows from theorem \ref{tn.3} by induction on $n$. $\bl$

\begin{thm}
\label{tn.5}
Let $0<p<1$ and let ${p_\#}$ be defined by {\em\rf{n.1}}.
Suppose that the block-matrix $\{\O_{kl}\}_{k.l\ge0}$ is block-triangular, i.e.,
$\O_{kl}=0$ for $k<l$. Then
$$
\sum_{k\ge0}\|\O_{kk}\|^{p_\#}_{\Mp}\le\|\O\|^{p_\#}_{\Mp}
$$
\end{thm}

\Pf Theorem \ref{tn.3} allows us to repeat word by word the second part of the 
proof of Theorem \ref{tn.2}. $\bl$

\begin{cor}
\label{tn.6}
Let $0<p<1$. Then for any positive integer $N$ we have
$$
\sum_{0\le l\le k\le l+N}\|\O_{kl}\|^{p_\#}_{\M_p}\le C(N,p)\|\O\|^{p_\#}_{\Mp}
$$ 
for some $C(N,p)$ depending only on $N$ and $p$.
\end{cor}

\Pf The result follows easily by induction on $N$. $\bl$

We obtain one more useful result on the Schur multiplier
norms of sums of matrices.

Recall that for $p\in(0,1)$, the number ${p_\#}$ is defined by \rf{n.1}.

\begin{thm}
\label{tn.7}
Let $0<p\le1$ and let $\{n_l\}_{l\ge0}$ be an increasing sequence 
in $\Z_+$ with $n_0=0$. For a matrix $A=\{a_{jk}\}_{j,k\ge0}$ consider
matrices $A_l=\left\{a_{jk}^{(l)}\right\}_{j,k\ge0}$, $l\ge1$, defined by
$$
a^{(l)}_{jk}=\left\{\begin{array}{ll}a_{jk},&n_{l-1}\le j<n_l,\\
0,&\mbox{otherwise}\,.
\end{array}\right.
$$
Then 
$$
\|A\|_{\Mp}\le\left(\sum\limits_{l\ge1}\|A_l\|_{\Mp}^{p_\fl}\right)^{1/p_\fl}.
$$
\end{thm}
 
\Pf Consider sequences $x=\{x_j\}_{j\ge0}$ and 
$y=\{y_k\}_{k\ge0}$ in $\ell^2$ of norm one.
We have to estimate $\|\{a_{jk}x_jy_k\}_{j,k\ge0}\|_{\bS_p}$.
Clearly,
$$
\|\{a_{jk}x_jy_k\}_{j,k\ge0}\|_{\bS_p}^p
\le\sum_{l\ge1}\left\|\left\{a_{jk}^{(l)}x_jy_k\right\}_{j,k\ge0}\right\|_{\bS_p}^p.
$$

Denote by $y^{(l)}=\left\{y^{(l)}_k\right\}_{k\ge0}$, $l\ge1$, 
the finitely supported sequence in $\ell^2$ defined by
$$
y^{(l)}_k=\left\{\begin{array}{ll}y_k,&n_{l-1}\le t<n_l,\\
0,&\mbox{otherwise}\,.
\end{array}\right.
$$

Clearly,
$$
\left\|\left\{a_{jk}^{(l)}x_jy_k\right\}_{j,k\ge0}\right\|_{\bS_p}
\le\left\|\left\{a_{jk}^{(l)}\right\}_{j,k\ge0}\right\|_{\Mp}
\left\|\left\{y^{(l)}\right\}_{l\ge0}\right\|_{\ell^2}
$$
Using H\"older's inequality with exponents $2/(2-p)$ and $2/p$, we obtain
\bey
\|\{a_{jk}x_jy_k\}_{j,k\ge0}\|_{\bS_p}^p
&\le&\left(\sum_{l\ge1}\left\|\left\{a_{jk}^{(l)}\right\}_{j,k\ge0}\right\|
_{\Mp}^{p_\fl}\right)^{1-p/2}
\left(\sum_{l\ge1}\|y_l\|_{\ell^2}^2\right)^{p/2}\\
&=&\left(\sum_{l\ge1}\left\|\left\{a_{jk}^{(l)}\right\}_{j,k\ge0}\right\|
_{\Mp}^{p_\fl}\right)^{1-p/2}
\eey
which proves the result. $\bl$

We complete this section by a description of the closure in $\Mp$ of the set of 
matrices with finitely many nonzero entries. We consider the following sets:
$$
c_{00}(\Z^2_+)\df\big\{\{a_{st}\}_{s,t\ge0}:~a_{st}=0\quad\mbox{except for finitely
many pairs}\quad (s,t)\big\}
$$
and
$$
c_0(\Z_+^2)\df\big\{\{a_{st}\}_{s,t\ge0}:~\lim_{s+t\to\be}a_{st}=0\big\}.
$$
We also denote by $\Mo$ the closure of $c_{00}(\Z^2_+)$ in the space $\Mp$. 
Clearly, \linebreak$\Mo\subset\Mp\cap c_0(\Z^2_+)$. It is obvious that 
${\frak M}_2^0=c_0(\Z^2_+)$. With any matrix \linebreak$A=\{a_{st}\}_{s,t\ge0}$ and 
positive integers $m$ and $n$ we associate the matrices 
\linebreak$A^{m,n}=\{a_{st}^{m,n}\}_{s,t\ge0}$ defined by
$$
a^{m,n}_{st}=\left\{\begin{array}{ll}
a_{st},&s\ge m,\ t\ge n,\\
0,&\mbox{otherwise}\,.
\end{array}\right.
$$

For $A\in\Mp$ it is easy to see that $A\in\Mo$ if and only if 
$\lim\limits_{m+n\to+\infty}A^{m,n}=0$ 
(in the space $\Mp$) (in fact it can easily be shown that if
$$
\lim_{n\to+\infty}A^{0,n}=0,\quad 
\lim_{n\to+\infty}A^{n,0}=0\quad\mbox{and}\quad\lim_{n\to+\infty}A^{n,n}=0
$$
in the space $\Mp$, then $A\in\Mo$).

It is clear now that if $\G_\psi\in\Mo$ for an analytic function $\psi$, then 
$$
\lim_{n\to\be}\|\G_{(S^*)^n\psi}\|_{\Mp}=0,
$$
where $S^*$ is backward shift , $(S^*f)(z)=(f-f(0))/z$.

Using the fact that the $\Mp$-norm is not equivalent to the $\ell^\infty$-norm for 
$p\not=2$, we may construct a sequence
of finite matrices $\{A_k\}_{k\ge0}$ such that $\|A_k\|_{\Mp}=1$ for all $k\ge0$ and 
$\lim\limits_{k\to+\infty}\|A_k\|_{\ell^\infty}=0$.
(Naturally, for a matrix $A=\{a_{st}\}_{0\le s,t\le j}$ by $\|A\|_{\ell^\be}$ we mean
$$
\sup\{|a_{st}|:~0\le s,\,t\le j\}).
$$
If we consider now the block-diagonal matrix
$$
\O=\left(\begin{array}{ccccc}
A_0&0&0&\cdots\\
0&A_1&0&\cdots\\
0&0&A_2&\cdots\\
\vdots&\vdots&\vdots&\ddots
\end{array}\right),
$$
then it is easy to see that $\O\notin\Mo$ and $\O\in \Mp\cap c_0(\Z^2_+)$ if $p\ge1$. 
Consequently, $\Mo\not=\Mp\cap c_0(\Z^2_+)$ for $p\in[1,2)\cup(2,+\infty)$.

Surprisingly enough for $p<1$ the situation is quite different.

\begin{thm}
\label{tn.8}
Let $0<p<1$. Then
$\Mo=\Mp\cap c_0(\Z^2_+)$.
\end{thm}

As before, with a matrix $\O=\{\o_{st}\}_{s,t\ge0}$ and an increasing sequence 
$\{n_k\}_{k\ge0}$ we associate the block-matrix $\{\O_{kl}\}_{k,l\ge0}$ defined
by \rf{n.2} with $m_k=n_k$.

\Pf Let $\O=\{\o_{st}\}_{s,t\ge0}\in c_0(\Z^2_+)\cap\Mp$. Suppose that $\O\notin\Mo$. 
We may assume that $\dist(\O,\Mo)>1$. Clearly, 
$\O-\O^{n,n}=(\O-\O^{n,0})+(\O^{n,0}-\O^{n,n})\in\Mo$. 
Hence, $\|\O^{n,n}\|_{\Mp}>1$ for any $n\ge0$. 
It follows that for each $m\in\Z_+$ there exists a positive integer $d$ such that
$$
\left\|\{\o_{st}\}_{m\le s,\,t<m+d}\right\|_{\Mp}\ge1.
$$
We can construct now by induction an increasing sequence $\{n_k\}_{k\ge0}$ such that
$n_0=0$ and $\|\O_{kk}\|_{\Mp}>1$ for all $k\in\Z_+$, and 
\bay
\label{n.5}
\|\O^{0,n_{k+1}}-\O^{n_k+1,n_{k+1}}\|_{\Mp}<2^{-k},\quad k\ge1
\ey
and
\bay
\label{n.6}
\|\O^{n_{k+1},0}-\O^{n_{k+1},n_k+1}\|_{\Mp}<2^{-k},\quad k\ge1.
\ey     

Consider now the matrix 
$\breve\O=\{\breve\o_{st}\}_{s,t\ge0}$ so that
$$
\breve\o_{st}=\left\{\begin{array}{ll}\o_{st},&n_k\le s<n_{k+1},~n_k\le t<n_{k+1},
\\0,&\mbox{otherwise}\,.\end{array}\right.
$$
It is easy to see that
$$
\breve\O=\O-\sum_{k\ge1}\left(\O^{0,n_{k+1}}-\O^{n_k+1,n_{k+1}}\right)
-\sum_{k\ge1}\left(\O^{n_{k+1},0}-\O^{n_{k+1},n_k+1}\right)
$$
and since by \rf{n.5} and \rf{n.6}
$$
\sum_{k\ge1}\left\|\O^{0,n_{k+1}}-\O^{n_k+1,n_{k+1}}\right\|^p_{\bS_p}+
\sum_{k\ge1}\left\|\O^{n_{k+1},0}-\O^{n_{k+1},n_k+1}\right\|^p_{\bS_p}<\be,
$$
it follows that $\breve\O\in\Mp$ which contradicts Theorem \ref{tn.2}. $\bl$

\

\setcounter{equation}{0}
\setcounter{section}{3}
\section{\bf Hankel--Schur Multipliers}

\

We study in this section Hankel--Schur multipliers of class $\bS_p$ for $0<p<1$.
\linebreak
We obtain various sharp necessary conditions and sufficient conditions. 
Then we obtain a 
characterization of the Hankel--Schur multipliers of $\bS_p$ for Hankel matrices
whose symbols have lacunary power series.

We start with the following
lemma that has been proved in [Pel3]. We give the proof here for completeness. 

\begin{lem}
\label{t3.1}
Let $\psi$ be an analytic polynomial of degree $m-1$. Suppose that
$0<p\le1$. Then
\bay
\label{3.1}
\|\G_\psi\|_{\Mp}\le(2m)^{1/p-1}\|\psi\|_p.
\ey
\end{lem}

Note that \rf{3.1} can be rewritten in the following form
$$
\|\G_\psi\|_{\Mp}\le(2m)^{1/{p_\#}}\|\psi\|_p,
$$
where {\it throughout this section ${p_\#}$ is always the number defined by} \rf{n.1}. 
Naturally, if $p=1$, we assume that ${p_\#}=\be$.

\Pf We have to prove that for any $B\in\bS_p$,
\beq
\label{3.0}
\|\G_\psi\star B\|_{\bS_p}\le(2m)^{1/{p_\#}}\|\psi\|_p\|B\|_{\bS_p}.
\end{equation}
In view of \rf{1.0}
it suffices to show that \rf{3.0} holds for operators $B$ of rank one.
Indeed, if $B$ is an arbitrary operator, then we can consider its Schmidt
expansion and apply \rf{3.1} to each term. 

Let $Bx=(x,\a)\b$, $x\in\ell^2$, where $\a=\{\a_n\}_{n\ge0}$, 
$\b=\{\b_n\}_{n\ge0}$.
Let $\z_j=e^{2\pi\text{i}j/(2m)}$, $0\le j\le 2m-1$. 
Define $\ell^2$ vectors $f_j$ and $g_j$, $0\le j\le2m-1$, by
$$
f_j(k)=\left\{\begin{array}{ll}
\overline{\psi(\z_j)}\overline \a_k\z_j^k,&0\le k<m,\\0,&k\ge m;\end{array}\right. 
\quad g_j(n)=\left\{\begin{array}{ll}
\b_n\bar\z_j^n,&0\le n<m,\\0,&n\ge m.\end{array}\right.
$$
We define the rank one operators $B_j$, $0\le j\le 2m-1$, by
$$
B_jx=(x,f_j)g_j, \quad x\in\ell^2.
$$

Let us show that
\beq
\label{3.2}
\G_\psi\star B=\frac{1}{2m}\sum_{j=0}^{2m-1}B_j.
\end{equation}
Indeed, if $k\ge m$ or $n\ge m$, then 
$(B_je_k,e_n)=((\G_\psi\star B)e_k,e_n)=0$.
On the other hand, if $n<m$ and $k<m$, then
$$
(B_je_k,e_n)=\psi(\z_j)\overline \a_k \b_n\bar\z_j^{n+k}
$$ 
and
$$
((\G_\psi\star B)e_k,e_n)=\overline \a_k\b_n\hat\psi(n+k).
$$
Identity \rf{3.2} follows now from the equality
$$
\hat\psi(d)=\frac{1}{2m}\sum_{j=0}^{2m-1}\psi(\z_j)\bar\z_j^{d},
\quad0\le d\le 2m-1,
$$
which is true for every polynomial $\psi$ of degree at most $2m-1$.

We have
$$
\|B_j\|_{\bS_p}\le\|\a\|_{\ell^2}\|\b\|_{\ell^2}|\psi(\z_j)|
=|\psi(\z_j)|\cdot\|B\|_{\bS_p}.
$$
It follows now from \rf{1.1} that
\beq
\label{3.3}
\|\G_\psi\star B\|^p_{\bS_p}\le
\frac{1}{(2m)^p}\|B\|_{\bS_p}^p\sum_{j=0}^{2m-1}|\psi(\z_j)|^p.
\end{equation}

For $\t\in\T$ we define the polynomial $\psi_\t$ by $\psi_\t(\z)=\psi(\t\z)$.
Let us show that 
$\|\G_{\psi_\t}\star B\|_{\bS_p}=\|\G_\psi\star B\|_{\bS_p}$ for
all $\t\in\T$. Indeed, this follows from the obvious equality
$$
\G_{\psi_\t}\star B=D_\t(\G_\psi\star B)D_\t,
$$
where $D_t$ is the unitary operator on $\ell^2$ defined by $D_\t e_n=\t^ne_n$
(here $\{e_n\}_{n\ge0}$ is the standard orthonormal basis of $\ell^2$).

To complete the proof, we integrate inequality \rf{3.3} in $\t$:
\bay
\|\G_\psi\star B\|^p_{\bS_p}
&=&\int_{\T}\|\G_{\psi_\t}\star B\|^p_{\bS_p}d\m(\t)
\le\frac{1}{(2m)^p}\|B\|_{\bS_p}^p
\sum_{j=0}^{2m-1}\int_{\T}|\psi_\t(\z_j)|^pd\m(\t)\nonumber\\
&=&\frac{2m}{(2m)^p}\|B\|_{\bS_p}^p\|\psi\|_p^p=
(2m)^{1-p}\|B\|_{\bS_p}^p\|\psi\|_p^p\nonumber
\ey
($\m$ stands for normalized Lebesgue measure on $\T$). $\bl$

{\bf Remark.} In the same way it can be proved that if $\psi$ is a {\it trigonometric}
polynomial of degree $m-1$ and $0<p<1$, then
\bay
\label{3.z}
\|\G_\psi\|_{\Mp}\le(4m)^{1/p-1}\|\psi\|_p
\ey
(see [Pel3]). In particular, if $F$ is any $C^\be$ function on $\R$ with compact support, 
then
\bay
\label{3.y}
\|\G_{F_{(m)}}\|_{\Mp}\le C(p,F).
\ey
Clearly, \rf{3.y} follows from \rf{3.z} and \rf{2.1}.

\begin{thm}
\label{t3.w}
Let $A\in\Mo$. Then 
\bay
\label{3.w}
\lim\limits_{m\to+\infty}A\star\G_{F_{(m)}}=F(0)A
\ey
in the space $\Mp$. 
\end{thm}

\Pf The result follows immediately from Theorem \ref{tn.8}, inequality 
\rf{3.y}, and the obvious fact that \rf{3.w} holds matrices $A$ with finitely many
nonzero entries. $\bl$

\begin{cor}
\label{t3v}
Let $A$ be a Hankel--Schur multiplier of $\bS_p$. Suppose that \linebreak$A\in\Mo$.
Then there exists a sequence $\{A_m\}_{m\ge1}$ in $c_{00}(\Z^2_+)$
such that $A_m$ is a Hankel matrix for every $m$ and
$\lim\limits_{m\to+\infty}A_m=A$ in the space $\Mp$. 
\end{cor}

\Pf We can consider a $C^\be$ function $F$ with compact support and such that $F(0)=1$.
Put now $A_m=A\st\G_{F_{(m)}}$. The result follows now from Theorem \ref{t3.w}. $\bl$

The following consequence of Lemma \ref{t3.1} has also been obtained in [Pel3]:
\beq
\label{3.4}
\|\G_\psi\|_{\bS_p}\le2^{1/p-1}m^{1/p}\|\psi\|_p
\end{equation}
for any analytic polynomial $\psi$ of degree $m-1$. Indeed, to prove \rf{3.4},
it is sufficient to write
$$
\|\G_\psi\|_{\bS_p}=\|\G_\psi\st Q^{(m-1)}\|_{\bS_p}
\le\|Q^{(m-1)}\|_{\bS_p}\|\G_\psi\|_{{\frak M}_p}
=m\|\G_\psi\|_{{\frak M}_p}
$$
and apply Lemma \ref{t3.1}. Recall that $Q^{(m-1)}$ is defined by \rf{1.4}.

We need the following consequence of the description of the Hankel operators of class
$\bS_p$ stated in \S 2.

\begin{thm}
\label{t3.2}
Let $0<p\le1$ and
let $\psi$ be a polynomial such that
$$
\hat\psi(j)=0\quad\mbox{unless}\quad2^{n-1}<j<2^{n+1}.
$$
Then
\beq
\label{3.5}
d_p2^{n/p}\|\psi\|_p\le\|\G_\psi\|_{\bS_p}\le D_p2^{n/p}\|\psi\|_p
\end{equation}
for some constants $d_p$ and $D_p$ not depending on $n$.
\end{thm}

\Pf The right inequality in \rf{3.5} is an immediate consequence of \rf{3.4}. 
To prove the left inequality in \rf{3.5}, we apply the description of Hankel operators
of class $\bS_p$ stated in \S 2. Clearly, under the hypotheses of the
lemma, $\psi*V_j=0$ for $j<n-1$ and for $j>n+1$. It follows that
\begin{eqnarray*}
\|\G_\psi\|_{\bS_p}^p&\ge&
\const2^n\left(\|\psi*V_{n-1}\|^p_p+\|\psi*V_n\|^p_p+\|\psi*V_{n+1}\|^p_p\right)\\
&\ge&\const2^n\|\psi*\left(V_{n-1}+V_n+V_{n+1}\right)\|^p_p
\end{eqnarray*}
(the second inequality is an immediate consequence of the triangle inequality
in $L^p$). Now it remains to observe that under the hypotheses of the theorem
$$
\psi*(V_{n-1}+V_n+V_{n+1})=\psi.\quad\bl
$$

{\bf Remark.} {\it It is easy to see that Theorem \ref{t3.2} can also be stated 
as follows. If $\psi$ is a polynomial such that
$$
\hat\psi(j)=0\quad\mbox{unless}\quad N\le j\le 2N.
$$
Then
$$
c_pN^{1/p}\|\psi\|_p\le\|\G_\psi\|_{\bS_p}\le C_pN^{1/p}\|\psi\|_p
$$
for some constants $c_p$ and $C_p$ not depending on $N$.}

Let us now obtain lower estimates for $\|\G_\psi\|_{\Mp}$. Consider
an infinitely differentiable function $r$ whose support is contained in
$[1/2,2]$. We define the polynomials $R_n$ by
\beq
\label{3.6}
R_n(z)=\sum_{k\ge0}r\left(\frac{k}{2^n}\right)z^k.
\end{equation}

\begin{thm}
\label{t3.3}
Let $0<p\le1$ and let $r$ be an infinitely differentiable function such that
$\supp r\subset[1/2,2]$. Then
$$
\|\G_\psi\|_{\Mp}\ge\const2^{n(1/p-1)}\|\psi*R_n\|_p=
\const2^{n/{p_\#}}\|\psi*R_n\|_p.
$$
where the $R_n$ are defined by {\em\rf{3.6}}.
\end{thm}

Note that the constants here and thereafter do not depend on $n$.

In particular, if we apply Theorem \ref{t3.3} to the polynomials $V_n$ defined
in \S 2, we obtain
\beq
\label{3.7}
\|\G_\psi\|_{\Mp}\ge\const2^{n/{p_\#}}\|\psi*V_n\|_p.
\end{equation}

\Pf We have
$$
\|\G_\psi\|_{\Mp}\ge\frac{\|\G_\psi\st\G_{R_n}\|_{\bS_p}}{\|\G_{R_n}\|_{\bS_p}}=
\frac{\|\G_{\psi*R_n}\|_{\bS_p}}{\|\G_{R_n}\|_{\bS_p}}.
$$
By Theorem \ref{t3.2},
$$
\frac{\|\G_{\psi*R_n}\|_{\bS_p}}{\|\G_{R_n}\|_{\bS_p}}\ge\const
\frac{\|\psi*R_n\|_p}{\|R_n\|_p}.
$$
Finally, by \rf{2.1},
$$
\frac{\|\psi*R_n\|_p}{\|R_n\|_p}\ge\const2^{n(1/p-1)}\|\psi*R_n\|_p
=\const2^{n/{p_\#}}\|\psi*R_n\|_p
$$
which completes the proof. $\bl$

Now we obtain a sharp estimate for the multiplier
norm $\|\G_\psi\|_{\Mp}$ for $\psi$ satisfying the hypotheses of Theorem
\ref{t3.2}.

\begin{thm}
\label{t3.4}
Let $0<p\le1$. Suppose that $\psi$ satisfies the hypotheses of Theorem
\ref{t3.2}. Then
\beq
\label{3.8}
\const2^{n/{p_\#}}\|\psi\|_p\le\|\G_\psi\|_{\Mp}\le\const2^{n/{p_\#}}\|\psi\|_p.
\end{equation}
\end{thm}

\Pf The right inequality in \rf{3.8} follows from Lemma \ref{t3.1}.
The left one follows easily from Theorem \ref{t3.3} but we give here a more
elementary argument. Consider the matrix $Q\df Q^{(2^{n+1})}$ defined in \rf{1.4}. 
Clearly, $\rank Q=1$ and $\|Q\|_{\bS_p}=2^{n+1}$. On 
the other hand it is obvious that $\G_\psi\st Q=\G_\psi$, and so
$$
\|\G_\psi\|_{\Mp}\ge\frac{\|\G_\psi\|_{\bS_p}}{\|\G_Q\|_{\bS_p}}=
2^{-n-1}\|\G_\psi\|_{\bS_p}\ge\const2^{-n}2^{n/p}\|\psi\|_p
$$
by \rf{3.5} which completes the proof. $\bl$

{\bf Remark.} {\it As in the case of $\bS_p$ norms it is easy to see that 
Theorem \ref{t3.4} can also be stated as follows. If $\psi$ is a polynomial such that
$$
\hat\psi(j)=0\quad\mbox{unless}\quad N\le j\le 2N,
$$
then
$$
c_pN^{1/{p_\#}}\|\psi\|_p\le\|\G_\psi\|_{\Mp}\le C_pN^{1/{p_\#}}\|\psi\|_p
$$
for some constants $c_p$ and $C_p$ not depending on $N$.}

Recall that for any analytic function $\psi$ in $\dd$ we have an expansion
\beq
\label{3.9}
\psi=\sum_{n\ge0}\psi*V_n.
\end{equation}
Theorem \ref{t3.4} gives a sharp estimate for the norm of $\G_{\psi*V_n}$
in $\Mp$. Roughly speaking, Theorem \ref{t3.4} resolves the problem of
describing the Hankel--Schur multipliers ``locally''.
We can obtain now the following ``global'' upper and lower estimates.

\begin{thm}
\label{t3.5}
Let $0<p\le1$ and let $\psi$ be a function analytic in $\dd$. Then
\beq
\label{3.10}
\const\cdot\sup_n2^{n(1-p)}\|\psi*V_n\|^p_p\le\|\G_\psi\|^p_{\Mp}\le
\const\sum_{n\in\Z_+}2^{n(1-p)}\|\psi*V_n\|^p_p,
\end{equation}
and so
\beq
\label{3.11}
\left(B_p^{1/{p_\#}}\right)_+
\subset\{\f:~\G_\f\in\Mp\}\subset\left(B_{p\be}^{1/{p_\#}}\right)_+.
\end{equation}
\end{thm}

\Pf The left inequality in \rf{3.10} is an immediate consequence of Theorem 
\ref{t3.3}, while the right one follows easily from \rf{3.9}, Lemma \ref{t3.1},
and inequality \rf{1.1}. Finally, \rf{3.11} follows immediately from
\rf{3.10} and the definitions of the Besov spaces given in \S 2. $\bl$

Note however, that the above necessary condition and sufficient condition 
do not lead to a description of Hankel--Schur multipliers of $\bS_p$ whose
symbols have lacunary Taylor series.

It follows from Theorem \ref{t3.5} that if $\{n_j\}_{j\ge0}$ is an 
{\it Hadamard lacunary sequence} of positive integers, i.e., 
\bay
\label{3.a}
\frac{n_{j+1}}{n_j}>\r>1\quad\mbox{for some}\quad\r\quad\mbox{and for}\quad j\ge0,
\ey
and $\psi=\sum\limits_{j\ge0}\l_jz^{n_j}$, then
\beq
\label{3.12}
\const\sup_{j\ge0}n_j^{1/{p_\#}}|\l_j|
\le\|\G_\psi\|_{\Mp}\le\const\sum_{j\ge0}n_j^{1/{p_\#}}|\l_j|.
\end{equation}
In other words,
$$
\ell^p\subset\left\{\{n_j^{1/{p_\#}}|\l_j|\}_{j\ge0}:
~\G_\psi\in\Mp\right\}\subset\ell^\be,\quad
\psi=\sum_{j\ge0}\l_jz^{n_j}.
$$

We are going to show in this section that 
$$
\left\{\{n_j^{1/{p_\#}}|\l_j|\}_{j\ge0}:
~\G_\psi\in\Mp\right\}=\ell^{p_\#},\quad
\psi=\sum_{j\ge0}\l_kz^{n_j}.
$$

First we obtain other sufficient conditions and necessary conditions. In
particular, we improve the sufficient condition given in Theorem \ref{t3.5}.

The following sufficient condition will be deduced from Theorem \ref{tn.1}.

\begin{thm}
\label{t3.h}
Let $0<p<1$ and let ${p_\#}$ be defined by {\em\rf{n.1}}. 
Suppose that $\f\in B_{p_\#}^{1/{p_\#}}$. Then $\G_\f\in\Mp$.
\end{thm}

\Pf This is an immediate consequence of Theorem \ref{tn.1} and the description of 
the Hankel operators of class $\bS_{p_\#}$ (see \S 2). $\bl$

Note that none of the sufficient conditions given by Theorems \ref{t3.5}
and \ref{t3.h} implies the other one.

We are going to obtain other sufficient condition later in this section. 
Meanwhile we obtain a necessary condition for $\G_\psi$ to be in $\Mp$ 
for a class of functions $\psi$ whose Taylor series have certain intervals of zeros.

\begin{thm}
\label{t3.i}
Suppose that $0<p<1$ and ${p_\#}$ is defined by {\em\rf{n.1}}. 
Let $\{\xi_k\}_{k\ge0}$ and $\{\eta_k\}_{k\ge0}$ be sequences of positive integers
such that
$$
\xi_k<\eta_k<\xi_{k+1},\quad \frac{\xi_{k+1}}{\eta_k}>d,
\quad\mbox{and}\quad \frac{\eta_k}{\xi_k}<D
$$
for some $d>1$ and $D>1$. Let $\psi$ be a function analytic in $\dd$ such that
$\hat\psi(j)=0$ for $\displaystyle j\not\in\bigcup_{k=0}^\be[\xi_k,\eta_k)$. If
$\G_\psi\in\Mp$, then $\psi\in B_{p\,{p_\#}}^{1/{p_\#}}$ and
$$
\|\psi\|_{B_{p\,{p_\#}}^{1/{p_\#}}}\le C(p,d,D)\|\G_\psi\|_{\Mp},
$$
where $C(p,d,D)$ can depend only on $p$, $d$, and $D$.
\end{thm}

We need the following lemma.

\begin{lem}
\label{t3.j}
Let $m\ge1$ and $N\ge0$ be integers and let $0<p<1$. Suppose that $f$ is a polynomial 
such that $\hat f(j)=0$ whenever $j<m$ or $j>m+N$. Consider the matrix 
$A=\{\hat f(k+l)\}_{0\le k<m+N,\,0\le l<m}$. Then
$$
\|f\|_p\le C\left(p,\frac{N}{m}\right)m^{-1/{p_\#}}\|A\|_{\Mp},
$$
where $C\left(p,\frac{N}{m}\right)$ depends only on $p$ and $\frac{N}{m}$.
\end{lem}

{\bf Proof of Lemma \ref{t3.j}.} Consider the infinite Hankel matrix
$$
\G_f=\{\hat f(k+l)\}_{k,l\in\Z}
$$ 
and for $\iota\in\Z_+$ we define its submatrix 
$A_\iota=\{\hat f(k+l)\}_{0\le k<m+N,\,\iota m\le l<(\iota+1)m}$. Clearly, 
$A_0=A$. It is easy to see that $\|A_\iota\|_{\Mp}\le\|A\|_{\Mp}$. Indeed, it is easy
to see that $A_\iota$ can be obtained from $A$ by deleting upper rows and adding zero
rows. By \rf{1.2}, 
\bay
\nonumber
\|\G_f\|_{\Mp}&\le&\left(\sum_{\iota\ge0}\|A_\iota\|^p_{\Mp}\right)^{1/p}
\le\left(\left(\frac{N}{m}+2\right)\|A\|_{\Mp}^p\right)^{1/p}\\
\label{3.i}
&\le&\left(\frac{N}{m}+2\right)^{1/p}\|A\|_{\Mp},
\ey
since at most $\frac{N}{m}+2$ terms in the above sum can be nonzero.
We have 
$$
f=\sum_{\kappa\ge0}f*V_\kappa
=\sum_{\{\kappa:\,2^{\kappa+1}>m\}}f*V_\kappa.
$$
Consider the matrix $Q^{(m+N)}$ defined by \rf{1.4}. We obtain
\bay
\label{3.j}
\|\G_f\|_{\Mp}\ge\frac{\|\G_f\|_{\bS_p}}{\|Q^{(m+N)}\|_{\bS_p}}
=\frac{\|\G_f\|_{\bS_p}}{m+N+1}.
\ey
By the description of the Hankel matrices of class $\bS_p$ (see \S 2), we have
\bey
\|\G_f\|^p_{\bS_p}&\ge&C(p)
\sum_{\{\kappa:\,2^{\kappa+1}>m\}}2^\kappa\|f*V_\kappa\|_p^p\\
&\ge&\frac{1}{2}C(p)m
\sum_{\{\kappa:\,2^{\kappa+1}>m\}}\|f*V_\kappa\|_p^p\\
&\ge&\frac{1}{2}C(p)m\|f\|_p^p,
\eey
where $C(p)$ is a constant that may depend only on $p$, and so by \rf{3.j},
$$
\|\G_f\|_{\Mp}\ge\frac{1}{2}C(p)\frac{m^{1/p}}{m+N+1}\|f\|_p.
$$
Hence,
\bey
\|f\|_p&\le&2(C(p))^{-1}\frac{m+N+1}{m^{1/p}}\|\G_f\|_{\Mp}\\&=&
2(C(p))^{-1}\left(m^{-1/{p_\#}}+\frac{N+1}{m}m^{-1/{p_\#}}\right)\|\G_f\|_{\Mp}.
\eey

The result follows now from \rf{3.i}. $\bl$

{\bf Proof of Theorem \ref{t3.i}.} Without loss of generality we may assume that
$$
\left(1-\frac{1}{d}\right)\xi_l>2\quad\mbox{and}\quad
\frac{(d-1)^2}{d}\xi_l>2\quad\mbox{for any}\quad l\in\Z_+.
$$
Put $m_l=\eta_{l-1}$ for $l\ge1$, $m_0\df0$, 
$m_l^\heartsuit=\left[\frac{1}{2}\left(1+\frac{1}{d}\right)\xi_l\right]$,
where $[x]$ denote the largest integer not exceeding $x$. By the above assumptions,
$m_l^\heartsuit>m_l$.

Consider the matrix $\O=\{\o_{st}\}_{s,t\ge0}$ defined by
$$
\o_{st}=\left\{\begin{array}{ll}\hat\psi(s+t),&m_l\le t<m_l^\heartsuit,
\quad \mbox{for some}\quad l\ge0,\\
0,&\mbox{otherwise}\,.
\end{array}\right.
$$
Clearly, $\|\O\|_{\Mp}\le\|\G_\psi\|_{\Mp}$.  

Put $n_k\df \xi_k-m_k^\heartsuit$ for $k\ge1$ 
and $n_0=0$. Let us show that $n_k<n_{k+1}$. Indeed,
$$
n_k\le\frac{\xi_k}{2}\left(1-\frac{1}{d}\right)+1<
\frac{\xi_k}{2}\left(1-\frac{1}{d}\right)d\le
\frac{\xi_{k+1}}{2}\left(1-\frac{1}{d}\right)\le n_{k+1},
$$
since $\frac{(d-1)^2}{2d}\xi_k>1$ and
by the hypotheses of the theorem, $d\xi_k\le\xi_{k+1}$.

We associate with the matrix $\O$ and the sequences
$\{m_l\}_{l\ge0}$ and $\{n_k\}_{k\ge0}$ the block-matrix $\{\O_{kl}\}_{k,l\ge0}$,
see \rf{n.2}. We are going to show that $\O_{kl}=0$ for $k<l$.

Indeed, suppose that $k<l$. Consider an entry $\o_{st}$ of $\O_{kl}$.
If $\o_{st}\neq0$, then the inequality $t<m_{k+1}$ 
implies $t<m_k^\heartsuit$. Then taking into account that $n_{k+1}\le n_l$, we obtain
$$
\eta_{l-1}=m_l\le s+t<m_l^\heartsuit+n_l\le\xi_l,
$$
and so $\l_{st}=\hat\psi(s+t)=0$.

Let us show that we can find a positive integer $N$ such that $n_{l+N}+m_l>\eta_l$
for all $l\in\Z_+$. Indeed, since
$$
n_{l+N}\ge\frac{1}{2}\left(1-\frac{1}{d}\right)\xi_{l+N}\quad\mbox{and}\quad
\frac{\xi_{l+N}}{\eta_l}\ge d^N,
$$
it is sufficient to choose $N$ so that 
$$
\frac{1}{2}d^N\left(1-\frac{1}{d}\right)>1.
$$
By Corollary \ref{tn.6},
\bay
\label{3.k}
\sum_{0\le l\le k\le l+N}\|\O_{kl}\|^{p_\#}_{\Mp}\le\const\|\O\|^{p_\#}_{\Mp}
\le\const\|\G_\psi\|_{\Mp}.
\ey

Given $k,l\in\Z_+$, we can associate with the matrix $\O_{kl}$ the infinite matrix
\linebreak$\breve\O_{kl}=\{\breve\o_{st}\}_{s,t\ge0}$ so that
$$
\breve\o_{st}=\left\{\begin{array}{ll}\o_{st},&n_k\le s<n_{k+1},~m_l\le t<m_{l+1},
\\0,&\mbox{otherwise}\,.\end{array}\right.
$$
Consider the matrices $\sum\limits_{k=l}^{l+N}\breve\O_{kl}$ and
$$
\O_l^\spadesuit\df\left\{\hat\psi(s+t)\right\}
_{\xi_l-m_l^\heartsuit\le s<\eta_l-m_l,\,m_l\le t<m_l^\heartsuit}.
$$
It is easy to see that 
$$
\|\O_l^\spadesuit\|_{\Mp}\le\left\|\sum_{k=l}^{l+N}\breve\O_{kl}\right\|_{\Mp},
$$
since $n_{l+N}>\eta_l-m_l$. Moreover,
$$
\left\|\sum_{k=l}^{l+N}\breve\O_{kl}\right\|_{\Mp}^{p_\#}\le C(p,N)
\sum_{k=l}^{l+N}\left\|\breve\O_{kl}\right\|^{p_\#}_{\Mp},
$$
where $C(p,N)$ is a constant that may depend only on $p$ and $N$.

Put
$$
f_l(z)=\sum_{j=\xi_l}^{\eta_l-1}\hat\psi(j)z^j.
$$
By Lemma \ref{t3.j}, we have
\bey
\sum_{l\ge0}\xi_l\|f_l\|^{p_\#}_p
&\le& c(p,d,D)\sum_{l\ge0}\|\O_l^\spadesuit\|^{p_\#}_{\Mp}\\
&\le&c(p,d,D)\sum_{l\ge0}\sum_{k=l}^{l+N}\|\O_{kl}\|^{p_\#}_{\Mp}
\le c(p,d,D)\|\G_\psi\|^{p_\#}_{Mp}
\eey
by \rf{3.k}, where $c(p,d,D)$ can depend only on $p$, $d$, and $D$. Note that
the constant in the conclusion of Lemma \ref{t3.j} can be chosen independent of $l$.
Indeed, this follows from the inequality
$$
\frac{\eta_l-m_l}{m_l^\heartsuit-m_l}\le\const
$$
which in turn follows from the fact that
$m_l^\heartsuit-m_l\ge\frac{1}{3}\left(1-\frac{1}{d}\right)\xi_l$ for 
sufficiently large $l$, and so $\eta_l-m_l\le \eta_l\le D\xi_l$.

To complete the proof, it is sufficient to use the definition of the Besov space 
$B_{p\,{p_\#}}^{1/{p_\#}}$ given in \S 2. $\bl$ 

{\bf Remark.} Note however that for an arbitrary analytic function $\psi$ the condition
$\G_\psi\in\Mp$ does not imply that $\psi\in B_{p\,{p_\#}}^{1/{p_\#}}$. Indeed, let 
$\psi=\sum\limits_{j\ge0}z^j$. Clearly, $\G_\psi\st A=A$ for any matrix $A$,
and so $\G_\psi\in\Mp$ for any $p>0$. On the other hand $\psi*V_n=V_n$ for any
$n\ge0$. By \rf{2.1}, 
$$
c2^{-n/{p_\#}}\le\|V_n\|_p\le C2^{-n/{p_\#}}\|V_n\|_p,\quad n\ge0
$$
for positive constants $c$ and $C$. Thus
$$
\sum_{n\ge0}2^n\|\psi*V_n\|^{p_\#}_p=\be,
$$
and so $\psi\not\in B_{p\,{p_\#}}^{1/{p_\#}}$ 
(see the definition of Besov spaces given in \S 2).
In fact among the Besov classes $B^{1/{p_\#}}_{pr}$, $0<r\le\be$, the function $\psi$ 
belongs only to $B^{1/{p_\#}}_{p\be}$.

\begin{cor}
\label{t3.u}
Under the hypotheses of Theorem \ref{t3.i} the condition $\G_\psi\in\Mp$ implies 
that $\G_\psi\in\Mo$
\end{cor}

\Pf In view of Theorem \ref{tn.8} it suffices to prove that 
$\lim\limits_{j\to+\infty}\hat\psi(j)=0$ for any $\psi\in B_{pp_\#}^{1/p_\#}$.

Let us first show that if $\psi\in\left(B_{p\be}^{1/p_\#}\right)_+$, then 
$\{\hat\psi(j)\}_{j\in\Z_+}\in\ell^\be$.
Suppose that $F$ is a $C^\be$ function on $\R$ such that $\supp F=[-1/2,1/2]$ and
$F(0)=1$. It follows easily from the definition of Besov spaces given in
\S 2 that 
$$
\|\psi*(z^nF_{(n)})\|_{H^p}^p\le \const\cdot n^{p-1},
$$
where $F_{(n)}$ is defined by \rf{2.z}. Clearly, the $n$th Fourier coefficient of
$\psi*(z^nF_{(n)})$ is $\hat\psi(n)$.
It remains to observe that $|\hat f(n)|^p\le (n+1)^{1-p}\|f\|_{H^p}^p$ for any $n\ge0$,
see [Pr], Ch. II, \S 11.

To complete the proof, we observe that $B_{pp_\#}^{1/p_\#}\subset B_{p\be}^{1/p_\#}$
and the polynomials are dense in $B_{pp_\#}^{1/p_\#}$ (see the definition of Besov classes
given in \S 2). $\bl$

Note also that to make the conclusion that the condition $\G_\psi\in\Mp$ implies
that $\psi\in B_{p\,{p_\#}}^{1/{p_\#}}$, 
it is not necessary that $\psi$ has many large intervals
on which $\hat\psi$ is zero. Indeed, the following result can easily be deduced
from Theorem \ref{t3.i}.

\begin{cor}
\label{t3.k}
Let $Q>1$ and $a<b<aQ$. Suppose that a function $\psi$ satisfies the conditions
$$
\hat\psi(j)=0,\quad\mbox{if}\quad j\quad\mbox{is even and}\quad
j\in\bigcup_{k=0}^\be\left[aQ^k,bQ^k\right]
$$
and
$$
\hat\psi(j)=0,\quad\mbox{if}\quad j\quad\mbox{is odd and}\quad
j\in\bigcup_{k=0}^\be\left[bQ^k,aQ^{k+1}\right].
$$
Then $\psi\in B_{p\,{p_\#}}^{1/{p_\#}}$ and 
$\|\psi\|_{B_{p\,{p_\#}}^{1/{p_\#}}}\le C(p)\|\G_\psi\|_{\Mp}$,
where $C(p)$ may depend only on $p$.
\end{cor}

\Pf Clearly,
$$
\|\G_{\psi(-z)}\|_{\Mp}=\|\G_\psi\|_{\Mp}. 
$$
Put $f(z)=(\psi(z)+\psi(-z))/2$ and $g=\psi-f$. Obviously, 
$\|\G_f\|_{\Mp}\le\const\|\G_\psi\|_{\Mp}$ and 
$\|\G_g\|_{\Mp}\le\const\|\G_\psi\|_{\Mp}$. It 
remains to observe that both $f$ and $g$ satisfy the hypotheses of Theorem 
\ref{t3.i}. $\bl$

The following result shows that Theorem \ref{t3.i} admits a converse under
a more restrictive condition on the sequences $\{\xi_k\}_{k\ge0}$ and 
$\{\eta_k\}_{k\ge0}$.

\begin{thm}
\label{t3.x}
Let $0<p<1$. Suppose that $\{\xi_k\}_{k\ge0}$ and 
$\{\eta_k\}_{k\ge0}$ are sequences of positive integers satisfying the condition
\bay
\label{3.x}
\sum_{k\ge1}\left(\frac{\eta_k-\xi_k+\eta_{k-1}}{\eta_k}\right)^{2/{p_\#}}<\be.
\ey
Suppose that $\psi$ is a function analytic in $\dd$ and such that $\hat\psi(j)=0$ 
for \linebreak$j\not\in\bigcup\limits_{k=0}^\be[\xi_k,\eta_k)$. 
If $\psi\in B_{p\,{p_\#}}^{1/{p_\#}}$, then $\G_\psi\in\Mp$.
\end{thm}

Note that if we could prove that the conclusion of Theorem \ref{t3.x} holds
for any sequences $\{\xi_k\}_{k\ge0}$ and $\{\eta_k\}_{k\ge0}$ satisfying the
hypotheses of Theorem \ref{t3.i}, it would follow from the definition of Besov
classes given in \S 2 that the condition $\psi\in B_{p\,{p_\#}}^{1/{p_\#}}$ 
is sufficient for $\G_\psi\in\Mp$ for and arbitrary $\psi$. 

{\bf Proof of Theorem \ref{t3.x}.} It is easy to see that \rf{3.x} implies that
$3\eta_{k-1}<\xi_k$ for sufficiently large $k$. We may assume that this is
true for all $k\ge1$. Put
$$
\psi_k=\sum_{j=\xi_k}^{\eta_k-1}\hat\psi(j)z^j.
$$

We define the matrices $A_k=\{\a^{(k)}_{st}\}_{s,t\ge0}$, 
$B_k=\{\b^{(k)}_{st}\}_{s,t\ge0}$, and 
$C_k=\{\g^{(k)}_{st}\}_{s,t\ge0}$ as follows:
$$
\a^{(k)}_{st}=\left\{\begin{array}{ll}\hat\psi(s+t),&\eta_{k-1}\le s,t<\eta_k,\\
0,&\mbox{otherwise}\,,
\end{array}\right.\quad k\ge0,
$$
$$
\b^{(k)}_{st}=\left\{\begin{array}{ll}
\hat\psi(s+t),&0\le s<\eta_{k-1},\,\eta_{k-1}\le t<\eta_k,\\
0,&\mbox{otherwise}\,,
\end{array}\right.\quad k\ge1,
$$
$$
\g^{(k)}_{st}=\left\{\begin{array}{ll}
\hat\psi(s+t),&\eta_{k-1}\le s<\eta_k,\,0\le t<\eta_{k-1},\\
0,&\mbox{otherwise}\,,
\end{array}\right.\quad k\ge1.
$$
Put $A=\sum\limits_{k\ge0}A_k$, $B=\sum\limits_{k\ge1}B_k$, 
$C=\sum\limits_{k\ge1}C_k$. Clearly, $\G_\psi=A+B+C$ and
$$
\|\G_\psi\|^p_{\Mp}\le\|A\|^p_{\Mp}+\|B\|^p_{\Mp}+\|C\|_{\Mp}^p.
$$
By Theorem \ref{tn.2}, 
$$
\|A\|_{\Mp}^{p_\#}=\sum_{k\ge0}\|A_k\|_{\Mp}^{p_\#}.
$$
It follows from Lemma \ref{t3.1} that 
$\|A_k\|_{\Mp}\le\const\eta_k^{1/{p_\#}}\|\psi_k\|_p$, and
so
$$
\|A\|_{\Mp}^{p_\#}\le\const\sum_{k\ge0}\eta_k\|\psi_k\|^{p_\#}_p
\le\const\|\psi\|^{p_\#}_{B_{p\,{p_\#}}^{1/{p_\#}}}
$$
(see the definition of Besov classes given in \S 2).

Again, by Lemma \ref{t3.1},
$$
\|B_k\|_{\Mp}=\|C_k\|_{\Mp}
\le\const(\eta_k-\xi_k+\eta_{k-1})^{1/{p_\#}}\|\psi_k\|_{L^p}.
$$

Recall that $p_\fl$ is defined by \rf{flat}. Applying Theorem \ref{tn.7}, we obtain
\bey 
\left\|\sum_{k\ge1}B_k\right\|_{\Mp}^{p_\fl}
&=&\left\|\sum_{k\ge1} C_k\right\|_{\Mp}^{p_\fl}
\le\sum_{k\ge1}\|B_k\|_{\Mp}^{p_\fl}\\
&\le&\sum_{k\ge1}(\eta_k-\xi_k+\eta_{k-1})^{{p_\fl}/{p_\#}}\|\psi_k\|_{L^p}^{p_\fl}\\
&\le&
\left(\sum_{k\ge1}\left(\frac{\eta_k-\xi_k+\eta_{k-1}}{\eta_k}\right)^{2/{p_\#}}\right)
^{1-{p_\fl}/{p_\#}}
\left(\sum_{k\ge1} \eta_k\|\psi_k\|_{L^p}^{p_\#}\right)^{{p_\fl}/{p_\#}}\\
&\le&\const
\left(\sum_{k\ge1}\left(\frac{\eta_k-\xi_k+\eta_{k-1}}{\eta_k}\right)^{2/{p_\#}}\right)
^{1-{p_\fl}/{p_\#}}
\|\psi\|_{B_{p\,{p_\#}}^{1/{p_\#}}}^p
\eey
by H\"{o}lder's inequality. $\bl$

We obtain now one more sufficient condition which allows us to improve the sufficient
condition in Theorem \ref{t3.5} as well as Theorem \ref{t3.h} in the case $p\le2/3$.
Recall that $p_\fl$ is defined by \rf{flat}.

\begin{thm}
\label{t3.y}
Let $0<p\le1$ and $p\le r\le\min\{1,p_\flat\}$.
Suppose that \linebreak$\psi\in B^{1/{p_\#}}_{r\,r_\fl}$.
Then $\G_\psi\in\Mp$.
\end{thm}

\Pf Consider the function $v$ in the definition of the polynomials $V_n$
involved in the definition of Besov spaces. We can represent $v$ as
$v=v^{(1)}+v^{(2)}+v^{(3)}$, where the functions $v^{(m)}$, $m=1,\,2,\,3$
are nonnegative infinitely differentiable and
$$
\supp v^{(1)}\subset[1/2,1],\quad\supp v^{(2)}\subset[3/4,3/2],\quad
\supp v^{(3)}\subset[1,2].
$$
Put
$$
V_n^{(m)}(z)= \sum_{k\ge1}v^{(m)}\left(\frac{k}{2^n}\right)z^k, \quad n\ge1,
\quad m=1,\,2,\,3.
$$
Then it can be seen from the definition of Besov spaces that the functions
$$
\psi^{(m)}=\sum_{k\ge1}\psi*V^{(m)}_k,\quad m=1,\,2,\,3,
$$
belong to $B^{1/{p_\#}}_{r\,r_\fl}$ and
$$
\psi=\psi*V_0+\psi^{(1)}+\psi^{(2)}+\psi^{(3)}.
$$
We prove that $\G_{\psi^{(2)}}\in\Mp$. The facts that $\G_{\psi^{(1)}}\in\Mp$
and $\G_{\psi^{(3)}}\in\Mp$ can be proved in the same way.

Moreover, we split $\psi^{(2)}$ in two functions
$$
\sum_{k\ge1}\psi*V^{(2)}_{2k}\quad\mbox{and}\quad\sum_{k\ge1}\psi*V^{(2)}_{2k-1}.
$$
We consider only the first function and it will be clear that for the second 
function the proof is the same.

In view of the said above we may suppose that
$\hat\psi(j)=0$ 
for $j\notin\bigcup\limits_{k=1}^\be(4^k,2\cdot4^k)$. 
Put $n_k=2\cdot4^k$ for $k\ge1$ and $n_0=0$.
Put
$$
\psi_k=\sum_{j=4^k}^{2\cdot4^k}\hat\psi(j)z^j.
$$
The condition that $\psi\in B^{1/{p_\#}}_{r\,r_\fl}$ means that
$$
\sum_{k\ge1}4^{kr_\fl/{p_\#}}\|\psi_k\|_r^{r_\fl}<\be.
$$

We define the matrices $A_k=\{\a^{(k)}_{st}\}_{s,t\ge0}$, 
$B_k=\{\b^{(k)}_{st}\}_{s,t\ge0}$, and 
$C_k=\{\g^{(k)}_{st}\}_{s,t\ge0}$ as follows:
$$
\a^{(k)}_{st}=\left\{\begin{array}{ll}\hat\psi(s+t),&n_{k-1}\le s,t<n_k,\\
0,&\mbox{otherwise}\,,
\end{array}\right.\quad k\ge1,
$$
$$
\b^{(k)}_{st}=\left\{\begin{array}{ll}
\hat\psi(s+t),&0\le s<n_{k-1},\,n_{k-1}\le t<n_k,\\
0,&\mbox{otherwise}\,,
\end{array}\right.\quad k\ge1,
$$
$$
\g^{(k)}_{st}=\left\{\begin{array}{ll}
\hat\psi(s+t),&n_{k-1}\le s<n_k,\,0\le t<n_{k-1},\\
0,&\mbox{otherwise}\,,
\end{array}\right.\quad k\ge1.
$$
Put $A=\sum\limits_{k\ge1}A_k$, $B=\sum\limits_{k\ge1}B_k$, 
$C=\sum\limits_{k\ge1}C_k$. Clearly, $\G_\psi=A+B+C$ and
$$
\|\G_\psi\|^p_{\Mp}\le\|A\|^p_{\Mp}+\|B\|^p_{\Mp}+\|C\|_{\Mp}^p.
$$
By Theorem \ref{tn.2},
$$
\|A\|_{\Mp}^{p_\#}=\sum_{k\ge1}\|A_k\|_{\Mp}^{p_\#}.
$$
It follows from Lemma \ref{t3.1} that
$$
\|A_k\|_{\Mp}\le\const n_k^{1/{p_\#}}\|\psi_k\|_p
\le\const n_k^{1/{p_\#}}\|\psi_k\|_r,
$$ 
and so
\bey
\|A\|_{\Mp}^{p_\#}&\le&\const\sum_{k\ge1}n_k\|\psi_k\|^{p_\#}_r\\
&\le&\const\left(\sum_{k\ge1}n_k^{r_\fl/{p_\#}}\|\psi_k\|
^{r_\fl}_r\right)^{{p_\#}/r_\fl}\\[.2cm]
&\le&\const\|\psi\|^{p_\#}_{B_{r\,r_\fl}^{1/{p_\#}}}
\eey
(see the definition of Besov classes given in \S 2)
because $\frac{r_\fl}{{p_\#}}\ge1$.

Clearly. $\|B\|_{\Mp}=\|C\|_{\Mp}$.
To estimate the norm of $B$, we consider sequences $x=\{x_s\}_{s\ge0}$ and 
$y=\{y_t\}_{t\ge0}$ in $\ell^2$ of norm one.

Let $B=\{\b_{st}\}_{s,t\ge0}$. We have to prove that
$$
\|\{\b_{st}x_sy_t\}_{s,t\ge0}\|_{\bS_p}\le
\const\|\psi\|_{B_{r\,r_\fl}^{1/{p_\#}}}.
$$
We may assume that $x_s\ge0$ for all $s\ge0$. 
Denote by $y^{(k)}=\{y^{(k)}_t\}_{t\ge0}$, $k\ge1$, the finitely supported sequence 
in $\ell^2$ defined by
$$
y^{(k)}_t=\left\{\begin{array}{ll}y_t,&\eta_{k-1}\le t<\eta_k,\\
0,&\mbox{otherwise}\,.
\end{array}\right.
$$
Put $\delta=2(\frac1p-\frac1r)$. Clearly, $0\le\d\le1$.
Note that 
$$
\sum\limits_{s=0}^{N-1}x_s^{2-2\d}
\le\left(\sum\limits_{s=0}^{N-1}x_s^2\right)^{1-\delta}N^\d.
$$

Let us first estimate $\|\{b_{st}x_s^{1-\d}y_t\}_{s,t\ge0}\|_{\bS_r}$.
Clearly,
$$
\b_{st}x_s^{1-\delta}y_t=\sum\limits_{k\ge1}\b^{(k)}_{st}x_s^{1-\d}y_t
$$ 
and
$$
\b^{(k)}_{st}x_s^{1-\delta}y_t=\b^{(k)}_{st}x_s^{1-\d}y^{(k)}_t,\quad k\ge1.
$$

We have
$$
\left\|\left\{\b_{st}x_s^{1-\delta}y_t\right\}_{s,t\ge0}\right\|_{\bS_r}^r\le
\sum_{k\ge1}
\left\|\left\{\b^{(k)}_{st}x_s^{1-\delta}y^{(k)}_t\right\}_{s,t\ge0}\right\|_{\bS_r}^r.
$$

It follows from Lemma \ref{t3.1} that
\bey
\left\|\left\{\b^{(k)}_{st}x_s^{1-\d}y^{(k)}_t\right\}_{s,t\ge0}\right\|_{\bS_r}
&\le&\const n_k^{\d/2}
\left\|\left\{\b_{st}^k\right\}\right\|_{{\frak M}_r}
\left\|y^{(k)}\right\|_{\ell^2}\\
&\le&\const n_k^{1/{p_\#}}\|\psi_k\|_r\left\|y^{(k)}\right\|_{\ell^2}.
\eey
Hence,
\bey
\left\|\left\{\b_{st}x_s^{1-\delta}y_t\right\}_{s,t\ge0}\right\|^r_{\bS_r}
&\le&\const\sum_{k\ge1}n_k^{r/{p_\#}}\|\psi_k\|^r_r
\left\|y^{(k)}\right\|_{\ell^2}^r\\
&=&\const\left(\sum_{k\ge1}
n_k^{r_\fl/{p_\#}}\|\psi_k\|^{r_\fl}_r\right)^{1-r/2}\\
&\le&\const\|\psi\|^r_{B^{1/{p_\#}}_{r\,r_\fl}}.
\eey
It remains to observe that $\left\|\{\b_{st}x_sy_t\}_{s,t\ge0}\right\|_{\bS_p}\le
\left\|\{\b_{st}x_s^{1-\delta}y_t\}_{s,t\ge0}\right\|_{\bS_r}$. $\bl$

\begin{cor}
\label{t3.s}
Let $0<p\le1$ and let $\psi\in\left(B^{1/{p_\#}}_{p\,p_\fl}\right)_+$. 
Then $\G_\psi\in\Mp$.
\end{cor}

\Pf The result immediately follows from Theorem \ref{t3.y} if we put $r=p$. $\bl$

Note that Corollary \ref{t3.s} improves the sufficient condition in Theorem \ref{t3.5}.

\begin{cor}
\label{t3.t}
Let $0<p\le2/3$ and let 
$\psi\in\left(B^{1/{p_\#}}_{p_\fl}\right)_+=
\left(B^{1/{p_\#}}_{p_\fl\,p_\fl}\right)_+$. 
Then $\G_\psi\in\Mp$.
\end{cor}

\Pf Put $r=p_\fl$. Since $p\le2/3$, $r\le1$. Now the result follows immediately
from Theorem \ref{t3.y}. $\bl$

Note that Corollary \ref{t3.t} improves Theorem \ref{t3.h} in the case $p\le2/3$.

We obtain another sufficient condition for $\G_\psi$ to be a Schur
multiplier of $\bS_p$. It follows from results of
Bo\v{z}ejko [Bo] that if $\psi\in B^0_{2\,\infty}$, then $\G_\psi\in{\frak M}_1$
(see also [Pis2] where a proof of Bo\v{z}ejko's result is given).
We use the method of [Pis2] to prove the following result.

\begin{thm}
\label{Boz}
Let $p\in(0,1]$ and let $\psi\in B^{1/p_\#}_{2\,p_\#}$. Then $\G_\psi\in\Mp$.
\end{thm}

\Pf The case $p=1$ is just the theorem of Bo\v{z}ejko. Consider the case $p<1$. Then $p_\flat<2$.
Let $a_n\df\hat\psi(n)$, $n\in\Z_+$. We define the matrices $\{\a^+_{jk}\}_{j,k\in\Z_+}$ and
$\{\a^-_{jk}\}_{j,k\in\Z_+}$ by
$$
\a^+_{jk}=\left\{\begin{array}{ll}a_{j+k},&j\ge k,\\
0,&\mbox{otherwise}\,,
\end{array}\right.
$$
and
$$
\a^-_{jk}=\left\{\begin{array}{ll}a_{j+k},&j< k,\\
0,&\mbox{otherwise}\,.
\end{array}\right.
$$
Let us show that $\{\a^+_{jk}\}_{j,k\in\Z_+}\in\Mp$. 
Consider sequences $x=\{x_j\}_{j\ge0}$ and \linebreak$y=\{y_k\}_{k\ge0}$ in $\ell^2$ of norm one.
We have to estimate $\left\|\{a^+_{jk}x_jy_k\}_{j,k\ge0}\right\|_{\bS_p}$. As in the proof
of Theorem \ref{tn.1} it is easy to show that
$$
\left\|\{a^+_{jk}x_jy_k\}_{j,k\ge0}\right\|_{\bS_p}
\le\left\|\{a^+_{jk}x_j\}_{j,k\ge0}\right\|_{\bS_{p_\flat}}\|y\|_{\ell^2}=
\left\|\{a^+_{jk}x_j\}_{j,k\ge0}\right\|_{\bS_{p_\flat}}.
$$
Thus it suffices to estimate
$\left\|\{a^+_{jk}x_j\}_{j,k\ge0}\right\|_{\bS_{p_\flat}}$. We have 
\bay
\label{1999}
\left\|\{a^+_{jk}x_j\}_{j,k\ge0}\right\|_{\bS_{p_\flat}}^{p_\flat}\le
\sum_{j=0}^{\be}\left(\sum_{k=0}^j|x_j|^2|a_{j+k}|^2\right)^{\frac{p_\flat}2}=
\sum_{j=0}^{\be}|x_j|^{p_\flat}\left(\sum_{k=0}^j|a_{j+k}|^2\right)^{\frac{p_\flat}2}.
\ey
Let $\lambda_n\df\left(\sum\limits_{k=2^n-1}^{2^{n+1}-2}|a_k|^2\right)^{1/2}$, $n\in\Z_+$.

Clearly, $\sum\limits_{k=0}^j|a_{j+k}|^2\le\lambda_n^2+\lambda_{n+1}^2$ if $2^n-1\le j<2^{n+1}-1$.
Thus \rf{1999} implies
\bay
\label{2000}
\nonumber
\left\|\{a^+_{jk}x_j\}_{j,k\ge0}\right\|_{\bS_{p_\flat}}^{p_\flat}
&\le&
\sum_{n=0}^\be\left(\sum_{j=2^n-1}^{2^{n+1}-2}|x_j|^{p_\fl}\right)(\l_n^{p_\fl}+\l_{n+1}^{p_\fl})\\[.2cm]
&\le&
\sum_{n=1}^\be\left(\sum_{j=0}^{2^n-1}|x_j|^{p_\fl}\right)(\l_{n-1}^{p_\fl}+\l_n^{p_\fl}).
\ey
Since $p<1$, we have by Hardy's inequality,  
$$
\sum_{n=1}^\be\left(\frac1n\sum_{j=0}^{n-1} |x_j|^{p_\flat}\right)^{2/p_\flat}\le
C(p)\sum_{j=0}^\be|x_j|^2=C(p).
$$ 
Consequently,
$$
\sum_{n=0}^\be\left(2^{-np_\flat/p_\#}\sum_{j=0}^{2^n-1}|x_j|^{p_\flat}\right)^
{2/p_\flat}\le C(p)\sum_{j=0}^\be|x_j|^2=C(p).
$$
It follows now from \rf{2000} and H\"older's inequality that
$$
\left\|\{a^+_{jk}x_j\}_{j,k\ge0}\right\|_{\bS_{p_\flat}}^{p_\flat}\le 
C(p)\left(\sum_{n=0}^\be2^n\l_n^{p_\#}\right)^{1-p_\flat/2}
\le C(p)\|\psi\|_{B^{1/p_\#}_{2\,p_\#}}^{p_\flat}.
$$
Hence, \mbox{$\{\a^+_{jk}\}_{j,k\in\Z_+}\in\Mp$}. In the same way
one can prove that $\{\a^-_{jk}\}_{j,k\in\Z_+}\in\Mp$.
This completes the proof. $\bl$

Note that Theorem \ref{Boz} improves Theorem \ref{t3.h} in the case $p>2/3$.

We have obtained various sharp necessary conditions and sufficient for $\G_\psi$ to
be a Schur multiplier of $\bS_p$. However, it remains unclear whether the
condition $\psi\in\B_{p\,{p_\#}}^{1/{p_\#}}$ itself is sufficient for $\G_\psi\in\Mp$.

Let us obtain now one more condition under which $\G_\psi\in\Mo$ whenever $\G_\psi\in\Mp$.

\begin{lem}
\label{taa} 
Suppose that $0<p<1$ and $\psi$ is a function analytic in $\dd$ such that $\G_\psi\in\Mp$. 
Let $m$, $n$ be nonnegative integers such that $m\le n$. Then
\bay
\label{qq}
|\hat\psi(n)|\le\left(\frac1{m+1}
\sum_{j=-m}^m\left(1-\frac{|j|}{m+1}\right)|\hat\psi(n+j)|^2\right)^{\frac{1-p}{2-p}}
\|\G_\psi\|_{\Mp}^{\frac p{2-p}}.
\ey
\end{lem}

\Pf It suffices to consider the case $m=n$ since
$\|\G_{(S^*)^{n-m}\psi}\|_{\Mp}\le\|\G_\psi\|_{\Mp}$.
Consider the matrix $Q^{(n)}$, see(1.5). Clearly, 
$$
\|\G_\psi\star Q^{(n)}\|_{{\bS}_2}^2=\sum\limits_{j=-n}^n(n+1-|j|)|\hat\psi(n+j)|^2.
$$ 
It is easy to see that 
$$
\|\G_\psi\star Q^{(n)}\|_{{\bS}_1}
\ge\left|\sum\limits_{j=0}^n(\G_\psi\star Q^{(n)}e_j,e_{n-j})\right|=(n+1)|\hat\psi(n)|.
$$
Consequently,
\bey
(n+1)|\hat\psi(n)|&\le&\|\G_\psi\star Q^{(n)}\|_{\bS_1}\le
\|\G_\psi\star Q^{(n)}\|_{\bS_2}^{\frac{2-2p}{2-p}}
\|\G_\psi\star Q^{(n)}\|_{\bS_p}^{\frac p{2-p}}\\
&\le&
\left(\sum_{j=-n}^n(n+1-j)|\hat\psi(n+j)|^2\right)^{\frac{1-p}{2-p}}
\left((n+1)\|\G_\psi\|_{\Mp}\right)^{\frac p{2-p}}
\eey
which implies \rf{qq}. $\bl$

For a subset $\L$ of $\Z_+$ we denote by $\r(\L)$ the density of $\L$, 
$$
\r(\L)\df\limsup_{n\to\be}\frac{\card([0,n]\cap\L)}{n+1}.
$$

\begin{thm}
\label{tab}
Suppose that $0<p<1$ and $\psi$ is a function analytic in $\dd$ such that $\G_\psi\in\Mp$. 
Let 
$$
\Lambda=\{k\in\Z_+:~\hat\psi(k)\not=0\}.
$$ 
Suppose that $\rho(\Lambda)=0$. Then $\G_\psi\in\Mp^0$. 
\end{thm}

\Pf It suffices to observe that Lemma \ref{taa} implies that 
$\lim\limits_{j\to\infty}\hat\psi(j)=0$ and apply Theorem \ref{tn.8}. $\bl$

{\bf Remark.} Note that for $p\ge1$ the situation is quite different. Indeed, if
\linebreak$\psi(z)=\sum\limits_{k\ge0}z^{2^k}$, then $\G_\psi\in\Mp\setminus\Mp^0$  
for all $p\ge1$ and $\rho(\Lambda)=0$. Note also that 
for $\psi(z)=\sum\limits_{k\ge0}z^{kN}$ we have
$\G_\psi\in\Mp\setminus\Mp^0$ for all $p>0$  and $\rho(\Lambda)=1/N$.

We proceed now to a characterization of the Hankel--Schur multipliers 
$\G_\psi$ of $\bS_p$, $0<p<1$, for Hadamard lacunary power series $\psi$.

\begin{thm}
\label{t3.l}
Let $0<p<1$. Suppose that
$\{n_j\}_{j\ge0}$ is an increasing sequence of integers whose terms can be represented
as a finite union of the terms of finitely many Hadamard lacunary sequences. Let
$\psi$ is a power series of the form
\linebreak$\psi=\sum\limits_{j\ge0}\l_kz^{n_j}$. 
Then $\G_\psi\in\Mp$ if and only if
\beq
\label{3.l}
\{n_j^{1/{p_\#}}|\l_j|\}_{j\ge0}\in\ell^{p_\#}.
\end{equation}
\end{thm}

Recall that we always assume that ${p_\#}$ is defined by \rf{n.1}.

\Pf Suppose that \rf{3.l} holds. It follows easily from the definition of
Besov spaces given in \S 2 that $\psi\in B_{p_\#}^{1/{p_\#}}$. 
It follows now from Theorem \ref{t3.h} that $\G_\psi\in\Mp$.

The converse follows from Theorem \ref{t3.i} and the following lemma.

\begin{lem}
\label{t3.m}
Under the hypotheses of Theorem \ref{t3.l} there exist sequences 
\linebreak$\{\xi_k\}_{k\ge0}$ and $\{\eta_k\}_{k\ge0}$ satisfying the hypotheses
of Theorem \ref{t3.i} such that
\bay
\label{3.m}
n_j\in\bigcup_{k\ge0}[\xi_k,\eta_k).
\ey
\end{lem}
 
Let us first complete the proof of Theorem \ref{t3.l}. By Lemma \ref{t3.m} and
Theorem \ref{t3.i}, $\psi\in B_{p\,{p_\#}}^{1/{p_\#}}$. 
It is easy to see from the definition of
Besov spaces given in \S 2 that for $\psi$ as in the hypotheses of Theorem \ref{t3.i}
this is equivalent to \rf{3.l}. $\bl$

{\bf Proof of Lemma \ref{t3.m}.} Without loss of generality we may assume that
$n_0>0$. Suppose that the set $\{n_k\}$, $k\ge0$, is the union
of the terms of $N$ Hadamard lacunary sequences each of which satisfies \rf{3.a} with 
the same number $\r>1$. Note that for any $j\ge0$ at least one of the numbers
$$
\frac{n_{j+1}}{n_j},\frac{n_{j+2}}{n_{j+1}},\cdots,\frac{n_{j+N+1}}{n_{j+N}}
$$
is greater than $\r^{1/N}$. Indeed, assume the contrary.
Then there exist $\iota_1$ and $\iota_2$ such that 
$$
1\le \iota_1<\iota_2\le N+1
$$ 
such that $n_{j+\iota_1}$ and
$n_{j+\iota_2}$ are terms of the same sequence satisfying \rf{3.a} which is impossible
since $\left(\r^{1/N}\right)^N=\r$. 

Denote by $\cal N$ the set of al $j\in\Z_+$ for which $\frac{n_{j+1}}{n_j}>\r^{1/N}$.
Clearly, $\cal N$ is an infinite set. We can enumerate the elements of $\cal N$ by
the increasing sequence $\{\eta_k-1\}_{k\ge0}$. If $n_j\in\cal N$ and $n_j=\eta_{k-1}$,
we define $\xi_{k+1}=n_{j+1}$. Finally, we put $\xi_0=1$.
It is easy to see that
$$
\frac{\eta_{k+1}-1}{\eta_k-1}>\r^{1/N},
$$
whence 
$$
\frac{\eta_{k+1}}{\eta_k}>\frac{\r^{1/N}+1}{2}.
$$
Clearly, the sequences $\{\xi_k\}_{k\ge0}$ and $\{\eta_k\}_{k\ge0}$ 
satisfy the hypotheses of Theorem \ref{t3.i} and \rf{3.m} holds. $\bl$

We conclude this section with a result on ``semi''-Hankel--Schur multipliers
of class $\bS_p$, $0<p\le1$, whose symbols have Hadamard lacunary power series.

For a function $\psi=\sum\limits_{j\ge0}\hat\psi(m)z^m$ we put 
$\G_\psi^-=\{\a_{jk}\}_{j,k\ge0}$, where
$$
\a_{jk}=\left\{\begin{array}{ll}\hat\psi(j+k),&j>k,\\0,&j\le k.\end{array}\right.
$$
We also put 
$$
\G^+_\psi\df\G_\psi-\G^-_\psi.
$$

\begin{thm}
\label{t3.n}
Let $0<p\le1$. Suppose that $\psi=\sum\limits_{j\ge0}\hat\psi(n_l)z^{n_l}$, where
$\{n_l\}_{l\ge0}$ is an increasing sequences of positive integers such that $n_0>0$
and $n_{l+1}\ge2n_l$, $l\in\Z_+$. Then
$$
\|\G^-_\psi\|_{\Mp}\le\|\G^-_\psi\|_{\ell^{p_\#}}
$$
where ${p_\#}=p/(1-p)$.
\end{thm}

Here for a matrix $B=\{\b_{jk}\}_{j,k\ge0}$ by $\|B\|_{\ell^{p_\#}}$ we mean
$\left(\sum\limits_{j,k\ge0}|\b_{jk}|^{p_\#}\right)^{1/{p_\#}}$ 
for \linebreak ${p_\#}<\be$ and
$\sup\limits_{j,k\ge0}|\b_{jk}|$ for ${p_\#}=\be$. 
Note that $\|\G^-_\psi\|_{\ell^{p_\#}}$ 
is equivalent to the norm of the sequence $\{n_l^{1/{p_\#}}\hat\psi(n_l)\}_{j\ge0}$ 
in $\ell^{p_\#}$.

We need a well-known inequality for $\bS_r$ norms of matrices. Let 
$A=\{a_{jk}\}_{j,k\ge0}$ be a matrix in $\bS_r$, $0<r\le2$. Then
\bay
\label{ineq}
\|A\|_{\bS_r}\le\|A\|_{\ell^r}=\left(\sum_{j,k\ge0}|a_{jk}|^r\right)^{1/r}.
\ey
Indeed, for $0<r\le1$ inequality \rf{ineq} is an immediate consequence of
\rf{1.0}, for $r=2$ we have equality, and for $1<r<2$ the result follows
by interpolation. 

{\bf Proof of Theorem \ref{t3.n}.} Consider a matrix $A=\{a_jb_k\}_{j,k\ge0}$ 
of rank one and assume that
$\|\{a_j\}_{j\ge0}\|_{\ell^2}=\|\{b_k\}_{k\ge0}\|_{\ell^2}=1$.
Consider the matrix $X=\{x_{jk}\}_{j,k\ge0}$ such that
$$
x_{jk}=\left\{\begin{array}{ll}a_j\l_l,&j+k=n_l\quad\mbox{for some}\quad k<j,\\
0,&\mbox{otherwise}\,,
\end{array}\right.
$$
where $\l_l\df\hat\psi(n_l)$. 

It is easy to see that the matrix $X$ is well-defined, since there can be at most
one such a number $k$. Indeed, if $j+k_1=n_{l_1}$ and $j+k_2=n_{l_2}$ and $l_1<l_2$, 
then $j>\frac{n_{l_2}}{2}\ge n_{l_1}$ and we get a contradiction.

Let $Y$ be the diagonal matrix with diagonal sequence $\{b_k\}_{k\ge0}$. It is easy
to see that 
$$
\G_\psi^-\st A=XY.
$$
Hence, 
\bey
\|\G_\psi^-\|_{\Mp}\le\|XY\|_{\bS_p}\le\|X\|_{\bS_{p_\fl}}\|Y\|_{\bS_2}=
\|X\|_{\bS_{p_\fl}}
\eey
(recall that $p_\fl$ is defined by \rf{flat}).

Let us now estimate $\|X\|_{\bS_{p_\fl}}$. Consider the set $\cal Q$ of
positive integers $j$ for which there exists $k<j$ such that $k+j=n_l$ for
some $l\in\Z_+$. For $j\in\cal Q$ we define $l(j)$ so that $l(j)<j$ and $l(j)+j=n_l$.
We have by \rf{ineq}
\bey
\|X\|_{\bS_{p_\fl}}^{p_\fl}&\le&
\sum\limits_{j\in{\cal Q}}|a_j|^{p_\fl}|\l_{l(j)}|^{p_\fl}\\
&\le&\left(\sum_{j\ge0}|a_j|^2\right)^{\frac{p}{2-p}}
\left(
\sum_{l\ge0}|\l_l|^{p_\#}\left[\frac{n_l+1}{2}\right]\right)^{\frac{2-2p}{2-p}}\\
&=&\|\G_\psi^-\|_{\ell^{p_\#}}^{p_\fl}
\eey
by H\"{o}lder's inequality.
Therefore $\|\G_\psi^-\|_{\Mp}\le\|\G^-_\psi\|_{\ell^{p_\#}}$. $\bl$

Note that if $n_{l+1}>2n_l$, $l\ge0$, then in a similar way one can prove that
\linebreak$\|\G_\psi^+\|_{\Mp}\le\|\G_\psi\|_{\ell^{p_\#}}$.

Theorem \ref{t3.n} straightforwardly implies a description of the Hankel--Schur multipliers 
of $\bS_1$ whose symbols have Hadamard lacunary Fourier series. This description is well known
(see [Bo], [Pis2], and [LP]).

\begin{cor}
\label{t3.o}
If $p=1$ and $\{n_l\}_{l\ge0}$ satisfies the hypotheses of Theorem \ref{t3.n}, then
$\|\G^-_\psi\|_{{\frak M}_1}=\|\G^-_\psi\|_{\ell^\be}$.
\end{cor}

\Pf By Theorem \ref{t3.n}, $\|\G^-_\psi\|_{{\frak M}_1}\le\|\G^-_\psi\|_{\ell^\be}$.
It remains to observe that obviously, 
$\|\G^-_\psi\|_{{\frak M}_1}\ge\|\G^-_\psi\|_{\ell^\be}$ for any function $\psi$.
$\bl$

\begin{cor}
\label{t3.q}
If $\psi$ and $\{n_l\}_{l\ge0}$ is a sequence whose terms are the union of the terms
finitely Hadamard lacunary sequences, then
$\G_\psi\in{\frak M}_1$ if and only if $\{\hat\psi(n_l)\}_{l\ge0}\in\ell^\be$.
\end{cor}

\Pf We have to prove that if $\{\hat\psi(n_l)\}_{l\ge0}\in\ell^\be$, then 
$\G_\psi\in{\frak M}_1$. The converse is obvious. 
Clearly, we can represents the terms of $\{n_l\}_{l\ge0}$ as the union of
the terms of finitely many sequences each of which satisfies the hypotheses
of Theorem \ref{t3.n}. It is easy to see that it is sufficient to
assume that the sequence $\{n_l\}_{l\ge0}$ itself satisfies the hypotheses
of Theorem \ref{t3.n}. By Theorem \ref{t3.n},
it follows that $\G_\psi^-\in{\frak M}_1$. Hence, $\G_\psi^+\in{\frak M}_1$,
and so $\G_\psi\in{\frak M}_1$. $\bl$

\

\setcounter{equation}{0}
\setcounter{section}{4}
\section{\bf Toeplitz--Schur multipliers}

\

In this section we use the results of the previous section to obtain
a description of the Toeplitz--Schur multipliers of $\bS_p$ for $0<p<1$.
Since $\Mp\subset{\frak M}_1$, it follows from the description of the 
Schur Toeplitz multipliers of $\bS_1$ (see \S 2) that if a Toeplitz matrix
$\{t_{j-k}\}_{j,k\ge0}$ is a Schur multiplier of $\bS_p$ for $0<p<1$, then
there exists a complex Borel measure $\mu$ on $\T$ such that $t_j=\hat\mu(j)$,
$j\in\Z$.

Let us introduce a class of measures. Let $0<p<1$. The class $\M_p$ is by
definition the space of discrete measures $\mu$ on $\T$ of the form
$$
\mu=\sum_j\a_j\d_{\t_j},\quad\a_j\in\C,~\t_j\in\T,\quad\mbox{the $\t_j$
are distinct},
$$
$$
\|\mu\|_{\M_p}\df\left(\sum_j|\a_j|^p\right)^{1/p}<\be.
$$

The main result of this section is the following theorem.

\begin{thm}
\label{t4.1}
Let $0<p<1$.
A Toeplitz matrix $\{t_{j-k}\}_{j,k\ge0}$ is a Schur multiplier of $\bS_p$
if and only if there exists $\mu\in\M_p$ such that $t_j=\hat\mu(j)$, $j\in\Z$.
Moreover, in this case
\beq
\label{4.*}
\|\{t_{j-k}\}_{j,k\ge0}\|_{\Mp}=\|\mu\|_{\M_p}.
\end{equation}
\end{thm}

\Pf The fact that the condition $\mu\in\M_p$ implies that
$\{\hat\mu(j-k)\}_{j,k\ge0}\in\Mp$ is very simple. Indeed, consider first
the case when $\mu=\d_\t$ for $\t\in\T$. We have by \rf{2.3},
$\hat\mu(j-k)=\bar\t^j\t^k$, and it is easy to see that 
$\|\{\hat\mu(j-k)\}_{j,k\ge0}\|_{\Mp}=1$.
Now if $\mu$ is an arbitrary measure in $\M_p$ of the form
$\mu=\sum_j\a_j\d_{\t_j}$, then it follows from \rf{1.1} that
\beq
\label{4.+}
\|\{\hat\mu(j-k)\}_{j,k\ge0}\|_{\Mp}^p
\le\sum_j|\a_j|^p\|\{\hat\d_{\t_j}(j-k)\}_{j,k\ge0}\|^p_{\Mp}
=\sum_j|\a_j|^p<\be.
\end{equation}

To prove the converse, we have to work harder. First we consider the class
of {\it Laurent matrices}. These are two-sided infinite matrices of the
form $\{t_{j-k}\}_{j,k\in\Z}$. 

We can consider two-sided Schur multipliers of
class $\bS_p$ that can be defined in the same way as one-sided:
$A=\{a_{jk}\}_{j,k\in\Z}$ is a Schur multiplier of $\bS_p$ if
$$
A\st B=\{a_{jk}b_{jk}\}_{j,k\in\Z}
$$
is a matrix of an operator of class $\bS_p$ 
on the two-sided sequence space $\ell^2(\Z)$ 
whenever $B=\{b_{jk}\}_{j,k\in\Z}\in\bS_p$.

\begin{lem}
\label{t4.2}
Let $\{t_j\}_{j\in\Z}$ be a two-sided sequence of complex numbers.
Then $\{t_{j-k}\}_{j,k\ge0}$ is a Schur multiplier of $\bS_p$
if and only if $\{t_{j-k}\}_{j,k\in\Z}$ is. Moreover,
$$
\|\{t_{j-k}\}_{j,k\ge0}\|_{\Mp}=\|\{t_{j-k}\}_{j,k\in\Z}\|_{\Mp}.
$$
\end{lem}

\Pf Suppose that the Laurent matrix $\{t_{j-k}\}_{j,k\in\Z}$ is a Schur multiplier of
$\bS_p$. Consider a matrix $B=\{b_{jk}\}_{j,k\in\Z}\in\bS_p$ such that
$b_{jk}=0$ whenever $j<0$ or $k<0$. It follows that the matrix
$\{t_{j-k}b_{jk}\}_{j,k\ge0}$ belongs to $\bS_p$, and so the Toeplitz matrix 
$\{t_{j-k}\}_{j,k\ge0}$ is a Schur multiplier of $\bS_p$ with
$\|\{t_{j-k}\}_{j,k\ge0}\|_{\Mp}\le\|\{t_{j-k}\}_{j,k\in\Z}\|_{\Mp}$.

Conversely, suppose that $\{t_{j-k}\}_{j,k\ge0}$ is a Schur multiplier of 
$\bS_p$. Consider the matrices $R_n=\{r^{(n)}_{jk}\}_{j,k\in\Z}$, $n\ge0$, defined by
$$
r^{(n)}_{jk}=\left\{\begin{array}{ll}
t_{j-k},&j\ge-n,~k\ge-n,\\0,&\mbox{otherwise}\,.\end{array}\right.
$$
It is obvious that $\|R_n\|_{\Mp}=\|\{t_{j-k}\}_{j,k\ge0}\|_{\Mp}$ for any
$n\in\Z_+$. It is also obvious that 
$$
R_n\st B\to\{t_{j-k}\}_{j,k\in\Z}\st B\quad\mbox{for any}\quad B\in\bS_p,
$$
and so $\{t_{j-k}\}_{j,k\in\Z}\in\Mp$ and 
$\|\{t_{j-k}\}_{j,k\in\Z}\|_{\Mp}=\|\{t_{j-k}\}_{j,k\ge0}\|_{\Mp}$. $\bl$

\begin{lem}
\label{t4.3}
Let $\{t_{j-k}\}_{j,k\ge0}$ be a Schur multiplier of $\bS_p$. Then for any
$m\in\Z$ the matrix $\{t_{j-k-m}\}_{j,k\ge0}$ is also a Schur multiplier of $\bS_p$
and
$$
\|\{t_{j-k-m}\}_{j,k\ge0}\|_{\Mp}=\|\{t_{j-k}\}_{j,k\ge0}\|_{\Mp}.
$$
\end{lem}

\Pf By Lemma \ref{t4.2}, the Laurent matrix $\{t_{j-k}\}_{j,k\in\Z}$ is a Schur 
multiplier of $\bS_p$. This is obvious that this is equivalent to the fact that
$\{t_{j-k-m}\}_{j,k\in\Z}$ is a Schur multiplier of $\bS_p$ and their
$\Mp$ quasi-norms are the same. Again by Lemma
\ref{t4.2}, this is equivalent to the fact that $\{t_{j-k-m}\}_{j,k\ge0}\in\Mp$.
$\bl$

\begin{lem}
\label{t4.4}
Let $\{t_{j-k}\}_{j,k\ge0}$ be a Schur multiplier of $\bS_p$. Then the Hankel matrix
$\{t_{j+k}\}_{j,k\ge0}$ is also a Schur multiplier of $\bS_p$ and
$$
\|\{t_{j+k}\}_{j,k\ge0}\|_{\Mp}\le\|\{t_{j-k}\}_{j,k\ge0}\|_{\Mp}.
$$
\end{lem}

\Pf By Lemma \ref{t4.2}, $\{t_{j-k}\}_{j,k\in\Z}\in\Mp$. Then its submatrix
$\{t_{j-k}\}_{j\ge0,k\le0}$ is also a Schur multiplier of $\bS_p$ and 
$$
\|\{t_{j-k}\}_{j\ge0,k\le0}\|_{\Mp}\le\|\{t_{j-k}\}_{j,k\ge0}\|_{\Mp}.
$$
Clearly, 
$$
\|\{t_{j+k}\}_{j,k\ge0}\|_{\Mp}=\|\{t_{j-k}\}_{j\ge0,k\le0}\|_{\Mp}
$$ 
which completes the proof. $\bl$

\begin{cor}
\label{t4.5}
Let $\{t_{j-k}\}_{j,k\ge0}$ be a Schur multiplier of $\bS_p$. Then for any $m\in\Z$
the Hankel matrix $\{t_{j+k-m}\}_{j,k\ge0}$ is also a Schur multiplier of $\bS_p$
and
$$
\sup_m\|\{t_{j+k-m}\}_{j,k\ge0}\|_{\Mp}\le\const.
$$
\end{cor}

\Pf The result follows immediately from Lemmas \ref{t4.3} and \ref{t4.4}. $\bl$

Consider now an infinitely differentiable even function $\o$ on $\R$ such that
\linebreak$0\le\o(s)\le1$, $s\in\R$, $\o(s)=1$ for $s\in[-1,1]$, and $\supp\o=[-2,2]$.
We define the trigonometric polynomials $\O_n$ by
\beq
\label{4.0}
\O_n(z)=\sum_{k\in\Z}\o\left(\frac{k}{2^n}\right)z^k,\quad n\ge1.
\end{equation}

Consider now a sequence $\{t_j\}_{j\in\Z}$ such that the Toeplitz matrix
$\{t_{j-k}\}_{j,k\ge0}$ is a Schur multiplier of $\bS_p$ for $0<p\le1$. As we
have observed above, there exists a complex Borel measure $\mu$ on $\T$
such that $t_j=\hat\mu(j)$, $j\in\Z$.

\begin{cor}
\label{t4.6}
Let $0<p\le1$ and let $\mu$ be a complex Borel measure on $\T$
such that $\{\hat\mu_{j-k}\}_{j,k\ge0}$ is a Schur multiplier of $\bS_p$.
Then 
\beq
\label{4.1}
\|\mu*\O_n\|_p\le\const2^{n(1-1/p)},\quad n\ge1.
\end{equation}
\end{cor}

\Pf By Theorem \ref{t3.3} and Lemma \ref{t4.4},
$$
\|\mu*z^{2^{n+2}}\O_n\|\le\const2^{n(1-1/p)},\quad n\ge1.
$$
Inequality \rf{4.1} follows now from Corollary \ref{t4.5}. $\bl$

We are going to deduce from \rf{4.1} the fact that
\beq
\label{4.2}
\|\mu*f\|_p\le\const\|f\|_p\quad\mbox{for any}\quad f\in L^1.
\end{equation}
To prove that $\mu\in\M_p$, it remains to use Oberlin's theorem [O] according to 
which \rf{4.2} implies that $\mu\in\M_p$. Note also that the proof
of the fact that \rf{4.1} implies that $\mu\in\M_p$ is essentially contained in
Theorem 2.3 of [A2], however we proceed in this paper via inequality \rf{4.2}
and Oberlin's theorem [O].

We need two more lemmas. For a function $g$ on $\R$ and $\e>0$ we put
$$
\f_\e(x)\df\frac{1}{\e}\f\left(\frac{x}{\e}\right),\quad x\in\R.
$$

\begin{lem}
\label{t4.7}
Let $0<p<1$ and $\d>0$. Suppose that 
$\f$ is a real function in $L^1(\R)\cap L^p(\R)$ and
$\int\limits_\R\f(x)dx\ne0$. Then for any function
$f\in L^1(\R)\cap L^p(\R)$ there exist $\a_j\in\C$, $\e_j\in(0,\delta)$, 
and $s_j\in\R$, $j\ge1$, such that
\beq
\label{4.3}
f(x)=\sum_{j\ge1}\a_j\f_{\e_j}(x-s_j),
\end{equation}
and
\beq
\label{4.4}
\sum_{j\ge1}|\a_j|^p\e_j^{1-p}\le C(p,\f)\|f\|_p^p
\quad\text{and}\quad
\sum_{j\ge1}|\a_j|\le C(\f)\|f\|_1.
\end{equation}
Moreover, one can choose $\e_j$ of the form $\e_j=2^{-n_j}$, $n_j\in\Z_+$.
\end{lem}

Note that \rf{4.4}  implies that the series in \rf{4.3} converges absolutely
in $L^1$ and $L^p$. Note also that implicitly this lemma is contained in [AK].

\Pf It suffices to prove that there exists $\s\in(0,1)$ such that for any \linebreak
$f\in L^1(\R)\cap L^p(\R)$ there exist sequences
$\a_j\in\C$, $0<\e_j<\delta$, and $s_j\in\R^d$ such that
$$
\int\limits_\R\left|f(x)-\sum_{j\ge1}\a_j\f_{\e_j}(x-s_j)\right|^rdx
\le(1-\s)\|f\|_r^r
$$
and
$$
\sum_{j\ge1}|\a_j|^r\e_j^{1-r}\le C(p,\f)\|f\|_r^r,\quad r=1,\,p.
$$
Clearly, it is sufficient to consider
the case $f=\chi_I$, where $I$ is an interval
in $\R$ ($\chi_I$ denotes the characteristic function of $I$).
Moreover, using translations and dilations, we see that it is sufficient
to prove that there exists a non-degenerate interval $I$ such that
the above inequality holds for $f=\chi_I$.
We may assume that $\int_\R\f(x)dx=1$. It follows easily from the Lebesgue
dominant convergence theorem that
\beq
\label{4.5}
\lim_{t\to0+}t^{-1}\int_\R(|1-t\f(x)|^r-1)dx=-r.
\end{equation}
Let us show that 
\beq
\label{4.6}
\|\chi_I-t\f\|_r^r<\|\chi_I\|_r^r,\quad r=1,\,p,
\end{equation}
if $I$ is sufficient large and $t$ is sufficient small.
We have
\begin{eqnarray*}
\|\chi_I-t\f\|_r^r-\|\chi_I\|_r^r&=&
\int_\R|\chi_I(x)-t\f(x)|^rdx-|I|\\
&=&\int_I(|1-t\f(x)|^r-1)dx+\int_{\R\setminus I}|t\f(x)|^rdx\\
&=&\int_\R(|1-t\f(x)|^r-1)dx\\
&+&\int_{\R\setminus I}(|t\f(x)|^r-|1-t\f(x)|^r+1)dx
\end{eqnarray*}
where $|I|$ stands for the length of $I$. It follows now from \rf{4.5} that
if $t$ is small enough, then $\int_\R(|1-t\f(x)|^r-1)dx<0$. Clearly, choosing $I$ 
sufficiently large we can make the modulus of 
$\int_{\R\setminus I}(|t\f(x)|^r-|1-t\f(x)|^r+1)dx$ as small as possible
which proves \rf{4.6}. 

It is clear from the proof that we can choose $\e_j$ to be of the form 
$\e_j=2^{-n_j}$, $n_j\in\Z_+$. $\bl$

In the following lemma we identify periodic functions on $[0,1]$ with functions
on $\T$ via the map $\xi(s)=e^{2\pi s\text{i}}$, $s\in[0,1]$. Let now 
$\f$ be the inverse Fourier transform of our function $\o$ defined after the proof of 
Corollary \ref{t4.5}, i.e.,
$$
\o(y)=\int_\R\f(x)e^{-2\pi\text{i}xy}dx,\quad y\in\R.
$$
Clearly, $\f$ satisfies the hypotheses of Lemma \ref{t4.7}. It is easy to see
that 
\beq
\label{4.8}
\O_n(x)=\sum_{j\in\Z}\f_{2^{-n}}(x+j),\quad x\in[0,1].
\end{equation}
Indeed, if $\O_n$ is the periodic function defined by \rf{4.8}, then its Fourier
coefficients are
$$
\hat\O_n(k)=\sum_{j\in\Z}2^n\int_0^1\f(2^nx+2^nj)e^{-2\pi\text{i}kx}dx
=\int_\R\f(s)e^{-2\pi\text{i}ks2^{-n}}ds=\o\left(\frac{k}{2^n}\right),
$$
i.e., the trigonometric polynomials defined by \rf{4.0}
coincide with the functions defined by \rf{4.8}.

\begin{lem}
\label{t4.8}
Let $0<p<1$ and let $\f$ and the $\O_n$ be as above.
Then any function $f$ on $\R$ with period 1 satisfying 
$f|[0,1]\in L^p[0,1]$ admits a representation
\beq
\label{4.9}
f(x)=\sum_{j\ge1}\a_j V_{n_j}(x-s_j)
\end{equation}
with some $\a_j\in\C$, $n_j\in\Z_+$, and $s_j\in\R$ such that
\beq
\label{4.10}
\sum_{j\ge1}|\a_j|^p 2^{n_j(p-1)}\le C(p,\f)\|f\|_p^p
\quad\text{and}\quad
\sum_{j\ge1}|\a_j|\le C(\f)\|f\|_{L^1}.
\end{equation}
\end{lem}

As in Lemma \ref{t4.7} the inequalities in \rf{4.10} imply that the series
\rf{4.9} converges absolutely in $L^1$ and $L^p$.

\Pf By Lemma \ref{t4.7}, we can represent the function
$f\chi_{[0,1]}$ in the form
$$
f\chi_{[0,1]}=
\sum_{j\ge1}\a_j\f_{2^{-n_j}}(x-s_j).
$$
Moreover, the series in \rf{4.10} converge. Then
$$
f(x)=\sum_{j\ge1}\a_j V_{n_j}(x-s_j),\quad x\in\R.\quad\bl
$$

Now we can complete the proof of Theorem \ref{t4.1}. 
Let $Tf=\mu*f$ for \linebreak$f\in L^1(\T)$.
It is clear that $T$ is a translation invariant operator. By Lemma \ref{t4.8},
$$
\|Tf\|^p_p\le\sum_{j\ge1}|\a_j|^p\|\mu*\O_{n_j}\|_p^p.
$$
By \rf{4.1},
$$
\|\mu*\O_{n_j}\|_p^p\le\const2^{n_j(p-1)},
$$
and so by \rf{4.10},
$$
\|Tf\|^p_p\le\const\sum_{j\ge1}|\a_j|^p2^{n_j(p-1)}\le\const\|f\|_p^p,\quad f\in L^1.
$$
Thus, the operator
$T:L^1(\T)\to L^1(\T)$ can be extended
to a translation-invariant operator on $L^p(\T)$.
Consequently, by Oberlin's theorem [O], $\mu\in\M_p$. 

Let us now prove \rf{4.*}. First of all, the inequality
$$
\|\{\hat\mu_{j-k}\}_{j,k\ge0}\|_{\Mp}\le\|\mu\|_{\M_p}
$$
has already been proved in \rf{4.+}. Let us establish the opposite inequality.
By Lemma \ref{t4.2}, it is sufficient to show that
\beq
\label{4.13}
\|\{\hat\mu_{j-k}\}_{j,k\in\Z}\|_{\Mp}\ge\|\mu\|_{\M_p}
\end{equation}
for any $\mu\in\M_p$. Since we have already proved that 
$\{\hat\mu_{j-k}\}_{j,k\ge0}\in\Mp$ if and only if $\mu\in\M_p$, it follows that
it is sufficient to prove \rf{4.13} for finite linear combinations of $\d$-measures.

Suppose that 
$$
\mu=\sum_{m=1}^N\a_m\d_{\t_m},\quad\a_m\in\C,\quad\t_m\in\T,\quad
\t_m\ne\t_l\quad\mbox{for}\quad m\ne l.
$$
Let $f$ be a function in $L^2$ such that $\|f\|_{L^2}=1$ and
$f_{\t_m}\perp f_{\t_l}$ for $m\ne l$, where $f_\t(\z)\df f(\t\z)$.
To construct such a function, it is sufficient to take any function
$f$ of norm 1 that is supported on a sufficiently small arc of $\T$.
Consider now the rank one matrix 
$$
A=\{\ov{\hat f(j)}\hat f(k)\}_{j,k\in\Z}.
$$
Clearly, $\|A\|_{\bS_p}=\|f\|_{L^2}^2=1$. We have
$$
\{\hat\mu_{j-k}\}_{j,k\in\Z}\st A=
\sum_{m=1}^N\a_m\bar\t_m^j\ov{\hat f(j)}\t_m^k\hat f(k)
=\sum_{m=1}^N\a_m\ov{\hat f_{\t_m}(j)}\hat f_{\t_m}(k),
$$
and so
$$
\|\{\hat\mu_{j-k}\}_{j,k\in\Z}\st A\|_{\bS_p}^p=
\sum_{m=1}^N|\a_m|^p=\|\mu\|_{\M_p}^p
$$
which completes the proof. $\bl$

Note that it is easy to see from the above argument that a for a sequence 
$\{t_j\}_{j\in\Z}$ the Toeplitz matrix $\{t_{j-k}\}_{j,k\ge0}$ is a Schur
multiplier of $\bS_p$, $0<p<1$,
if and only if for any $m\in\Z_+$ the Hankel matrix $\{t_{j+k-m}\}_{j,k\ge0}$
is a Schur multiplier of $\bS_p$ and
$$
\sup_{m\in\Z}\|\{t_{j+k-m}\}_{j,k\ge0}\|_{\Mp}\le\const.
$$
However, it can be shown easily that such a relation between Toeplitz and 
Hankel--Schur multipliers holds in a more general situation. 
In particular, this is also true
for Toeplitz and Hankel--Schur multipliers of $\bS_p$ with any $p\in(0,\be)$.

{\bf Remark.} We can obtain similar results for Toeplitz (Wiener--Hopf) 
multipliers of class $\bS_p$ of operator on $L^2(\R_+)$. Namely, we can ask the 
question of when the function $(x,y)\mapsto g(x-y)$, $x,\,y\ge0$, is a Schur 
multiplier of $\bS_p$. The techniques described in this paper allows us to prove that
this condition is equivalent to the fact that $g$ is a Fourier transform
of a complex discrete measure $\mu$ on $\R$ such that
$$
\mu=\sum_j\a_j\d_{x_j},\quad\mbox{the $x_j$ are distinct, and}
\quad\sum_j|\a_j|^p<\be.
$$

\

\setcounter{equation}{0}
\setcounter{section}{5}
\section{\bf More results on Hankel--Schur multipliers}

\

In this section we study Hankel--Schur multipliers of the form $\G_\mu$, where $\mu$ is a 
complex measure on $\T$. 
Recall that the space ${\cal M}_p$ has been introduced in \S 5.

If $\mu$ is a complex measure on the unit circle $\T$, we identify $\pp_+\mu$ and $\pp_-\mu$
with functions $\mu_+$ and $\mu_-$ defined by
$$
\mu_+(z)\df\int\limits_{\T}\frac{d\mu(\zeta)}{(1-\overline\zeta z)},\quad|z|<1
$$
and 
$$
\mu_-(z)\df-\int\limits_{\T}\frac{d\mu(\zeta)}{(1-\overline\zeta z)},\quad|z|>1.
$$

We need the following well known facts about Besov classes (we refer the reader to [Pee] for
more information about Besov spaces). Let $f\in H^p$, $0<p<\be$,
and let $0<s<1$. Then 
\bay
\label{6.1}
f\in B^s_{p\,\infty}\quad\Longleftrightarrow\quad
\int\limits_{\T}|f(\z\t)-f(\t)|^pd\mu(\t)\le\const|\z-1|^{sp},\quad\z\in\T,
\ey
and
\bay
\label{6.2}
f\in B^s_{p\,\infty}\quad\Longleftrightarrow\quad
\int\limits_{\T}|f(\z)-f(r\z)|^pd\m(\z)\le\const (1-r)^{sp}.
\ey

\begin{thm}
\label{t6.1}
Let $\mu$ be a singular measure on the unit circle $\T$. Suppose that 
$\mu_+\in B^s_{p\,\infty}$ , where
$p,s\in(0,1)$. Then $\mu_-\in B^s_{p\,\infty}$.
\end{thm}

\Pf By \rf{6.1},
$$
\int\limits_{\T}|\mu_+(\z\t)-\mu_+(\t)|^pd\mu(\t)\le\const|\z-1|^{sp},\quad\z\in\T.
$$
Since $\mu$ is singular, the boundary values of $\mu_+$ and $-\mu_-$ coincides almost 
everywhere on $\T$, and so.
$$
\int\limits_{\T}|\mu_-(\z\t)-\mu_-(\t)|^pd\mu(\t)\le\const|\z-1|^{sp},\quad\z\in\T,
$$
which implies by \rf{6.1} that $\mu_-\in B^s_{p\,\be}$. $\bl$

\begin{cor}
\label{t6.2}
Let $\mu$ be a singular measure on the unit circle $\T$. Suppose that 
$\mu_+\in B^{1/p_\#}_{p\,\infty}$ and $1/2<p<1$.
Then $\mu\in\M_p$.
\end{cor}

\Pf Denote by $u$ the Poisson integral of the measure $\mu$. Clearly,
$$u(z)=\mu_+(z)+\mu_-(\frac1{\overline z}),\quad z\in\dd.
$$
Let us estimate the integral $\int\limits_{\T}|u(r\zeta)|^p d\m(\zeta)$. 
Denote by $h$ the boundary values of $\mu_+$. They
coincide almost everywhere with the boundary values of $-\mu_-$ since $\mu$ is singular.
Note that by \rf{6.2},
$$
\int\limits_{\T}|\mu_+(r\zeta)-h(r\zeta)|^p d\m(\z)\le\const (1-r)^{1-p}.
$$
Similarly,
$$
\int\limits_{\T}|\mu_-(r^{-1}\zeta)+h(r^{-1}\zeta)|^p d\m(\z)\le\const (1-r)^{1-p}.
$$ 
Hence,
\bey
\int\limits_{\T}|u(r\zeta)|^p d\m(\zeta)
&\le&\int\limits_{\T}|\mu_+(r\zeta)-h(r\zeta)|^p d\m(\z)\\
&+&\int\limits_{\T}|\mu_-(r^{-1}\zeta)+h(r^{-1}\zeta)|^p d\m(\zeta)\\
&\le&\const (1-r)^{1-p}.
\eey
Consequently, by Theorem 4.3 of [A2], we have $\mu\in\M_p$. $\bl$

\begin{cor} 
\label{t6.3}
Let $\mu$ be a singular measure on $\T$. Suppose that $\G_\mu\in\Mp$ for some $p<1$.
Then $\mu$ is discrete.
\end{cor}

\Pf Clearly, $\G_\mu\in{\frak M}_q$ for all $q\ge p$. Consequently,  it suffices to 
consider the case $1/2<p<1$. By Theorem \ref{t3.5}, $\mu_+\in B^{1/p_\#}_{p\,\infty}$. 
It remains to apply Corollary \ref{t6.2}. $\bl$

It is clear from the proof of Corollary \ref{t6.3} that 
$\mu\in\M_p$ if $1/2<p<1$. We will see that the same is also true for $p\le1/2$.

Denote by $S^*$ backward shift, $(S^*\f)(z)\df (\f(z)-\f(0))/z$. It is easy to see 
that $\|\G_{S^*\psi}\|_{\Mp}\le\|\G_\psi\|_{\Mp}$ for any
$p>0$. Moreover, it is easy to see that if 
$f(z)=\sum\limits_{j\ge0}\hat\psi(jN+s)z^j$ with $N,s\in\Z_+$, then
\bay
\label{6.3}
\|\G_f\|_{\Mp}\le\|\G_\psi\|_{\Mp},\quad p>0.
\ey
Indeed, to prove \rf{6.3}, it is suffices to observe that the matrix $\G_f$ is a 
submatrix of $\G_\psi$.

\begin{thm}
\label{t6.4}
Let $\psi$ and $g$ be functions analytic in the unit disk $\dd$ and let
\linebreak$0<p<\be$. Suppose that $\{n_j\}_{j\ge0}$ is an increasing sequence
of nonnegative integers, $\G_\psi\in\Mp$, and 
$$
\lim\limits_{j\to+\infty}\left((S^*)^{n_j}\psi\right)(z)=g(z),\quad z\in\dd.
$$  
Then $\G_g\in\Mp$ and $\|\G_g\|_{\Mp}\le\|\G_\psi\|_{\Mp}$. 
Moreover, there exists a Schur multiplier $\{t_{k+l}\}_{k,l\in\Z}$ of $\bS_p$ such that
$t_k=\hat g(k)$ for all $k\ge0$ and 
$\|\{t_{k+l}\}_{k,l\in\Z}\|_{\Mp}\le\|\G_\psi\|_{\Mp}$.
\end{thm}

\Pf First we can put $\hat\psi(k)=0$ for $k<0$.
Clearly, $\lim\limits_{j\to+\infty}\widehat\psi(k+n_j)=\hat g(k)$ for any 
$k\ge0$. Passing to a subsequence (if necessary),
we may assume that sequence $\{\widehat\psi(k+n_j)\}_{j\ge0}$ converges for any $k\in\Z$. 
Put $t_k\df\lim\limits_{j\to+\infty}\widehat\psi(k+n_j)$
for $k\in\Z$. We have to prove that $\{t_{k+l}\}_{k,l\in\Z}\in\Mp$ and 
$\|\{t_{k+l}\}_{k,l\in\Z}\|_{\Mp}\le\|\G_\psi\|_{\Mp}$. 
Let $a=\{a_k\}_{k\in\Z}$ and $b=\{b_l\}_{l\in\Z}$ be two 
sequences in $\ell^2(\Z)$ of norm one such that 
$a_k=b_l=0$ if $|k|,|l|>N$ for some $N\in\Z_+$.
It is easy to see that
\bey
\left\|\left\{\hat\psi(k+l+n_j)a_kb_l\right\}_{k,l\in\Z}\right\|_{\bS_p}&=&
\left\|\left\{\hat\psi(k+l)a_{k+[n_j/2]} b_{l+n_j-[n_j/2]}\right\}_{k,l\ge0}\right\|_{\bS_p}\\
[.2cm]&\le&\|\G_\psi\|_{\Mp},
\eey 
if $n_j>2N$. Making $j\to\infty$, 
we obtain $\|\{t_{k+l}a_kb_l\}_{k,l\in\Z}\|_{\bS_p}\le\|\G_\psi\|_{\Mp}$. $\bl$

In the same way we may prove the following result.

\begin{thm}
\label{t6.5}
Let $g$ and $\psi_j$, $j\in\Z_+$, be  
functions analytic in the unit disk $\dd$ and let $0<p<\be$.
Suppose that $M=\sup\limits_{j\ge0}\|\G_{\psi_j}\|_{\Mp}<\be$ and
for an increasing sequence $\{n_j\}_{j\ge0}$ of positive integers
$\lim\limits_{j\to+\infty}\left((S^*)^{n_j}\psi_j\right)(z)=g(z)$ for all $z\in\dd$.
Then $\G_g\in\Mp$ and $\|\G_g\|_{\Mp}\le M$. Moreover, there exists 
a Schur multiplier $\{t_{k+l}\}_{k,l\in\Z}$ of $\bS_p$ such that
$t_k=\hat g(k)$ for all $k\ge0$ and $\|\{t_{k+l}\}_{k,l\in\Z}\|_{\Mp}\le M$.
\end{thm}

Denote by $\M(\T)$ the space all complex Borel measures on $\T$. As usual, we identify 
$\M(\T)$ with the space $(C(\T))^*$, where $C(\T)$ is the space of all continuous
functions on $\T$. Recall that a measure $\mu\in\M(\T)$ is called {\it continuous}
if $\mu\{\zeta\}=0$ for any $\zeta\in\T$.

\begin{lem} 
\label{t6.6}
Let $\mu$ be a continuous measure and $\nu$ a discrete measure in $\M(\T)$.
Then there exists an increasing sequence $\{n_j\}_{j\ge0}$ in $\Z_+$ such that 
$\lim\limits_{j\to\infty}\hat\mu(k+n_j)=0$ for any $k\in\Z$
and $\lim\limits_{j\to\infty}z^{n_j}=1$ for $|\nu|$-almost all $z\in\T$. 
\end{lem}

\Pf Assume first that the support of $\nu$ is finite.  
By Wiener's theorem [W], $\mu$ is a continuous measure if and only if 
$$
\lim_{N\to+\be}\frac1{N+1}\sum_{k=0}^N|\hat\mu(k)|^2=0.
$$
Consequently, there exists a sequence $\{m_j\}_{j\ge0}$ in $\Z_+$ such that $m_j>2^j$ and 
$$
\lim\limits_{j\to\infty}\sup\limits_{m_j-2j\le k\le m_j+2j}|\hat\mu(k)|=0
$$
Let us prove now that there exists a required sequence $\{n_j\}_{j\ge0}$ such that 
\bay
\label{6.4}
m_j-j\le n_j\le m_j+j,\quad j\in\Z_+.
\ey
Indeed, suppose that $\supp\nu=\{\zeta_1, \zeta_2, ...,  \zeta_N\}$. 
Let $X$ be the closure of the set $\{(\zeta_1^n,\zeta_2^n,...,\zeta_N^n)\}_{n\in \Z}$. 
Clearly, $X$ is a subgroup of $\T^N$. Set 
$$
X_n\df\{(\zeta_1^k,\zeta_2^k,...,\zeta_N^k):~-n\le k\le n\}.
$$ 
Denote by $\e_n$  the infimum of the set of positive $\e$ such that the
$\e$-neighborhood of $X_n$ contains $X$. Clearly, 
$\lim\limits_{n\to\infty}\e_n=0$. Now it is easy to choose a sequence $\{n_j\}_{j\ge0}$
satisfying \rf{6.4} and such that 
$\lim\limits_{j\to\infty}z^{n_j}=1$ for $|\nu|$-almost all $z\in\T$. 

If $\nu=\sum_{j=0}^\be c_j\d_{\t_j}$ is an arbitrary discrete measure, we can take 
the finitely supported measures $\nu_m=\sum_{j=0}^m c_j\d_{\t_j}$, apply the above
reasoning to the $\nu_m$, use a diagonal process. $\bl$

\begin{thm} 
\label{t6.7}
Let $\mu\in\M(\T)$ and $0<p<1$. Suppose that $\G_{\mu_+}\in\Mp$. 
Then $\sum\limits_{\zeta\in\T}|\mu\{\zeta\}|^p\le\|\G_{\mu_+}\|_{\Mp}^p$.
\end{thm}

\Pf Denote by $\nu$ the discrete part of measure $\mu$. By Lemma \ref{t6.6}, there exists an 
increasing sequence $\{n_j\}_{j\ge0}$ in $\Z_+$
such that $\lim\limits_{j\to\infty}\hat\mu(k+n_j)=\hat\nu(k)$. Theorem \ref{t6.4} implies that 
$\{\hat\nu_{k+l}\}_{k,l\in\Z}\in\Mp$ and
$\|\{\hat\nu_{k+l}\}_{k,l\in\Z}\|_{\Mp}\le\|\G_{\mu}\|_{\Mp}$.
It remains to apply Theorem \ref{t4.1}. $\bl$

\begin{thm} 
\label{t6.8}
Let $\mu$ be a singular measure on $\T$. Suppose that $\G_\mu\in\Mp$, where $0<p<1$.
Then $\mu\in\M_p$.
\end{thm}

\Pf It suffices to note that by Corollary \ref{t6.3} the measure $\mu$ is discrete. $\bl$ 

Finally, we consider the problem of whether all Hankel--Schur multipliers of $\bS_p$ for $0<p<1$
are of the form $\G_\mu$ for $\mu\in\M(\T)$. We show that this is not the case for $p>2/3$.
For $p<2/3$ the question remains open. 

We denote by $\M_+$ the space of Cauchy integrals of measures:
$$
\M_+\df\left\{\f:~\f(\z)=\int\limits_\dd\frac{d\mu(\t)}{1-\bar{\t}\z},~\z\in\dd,~
\mu\in\M(\T)\right\}
$$
and consider the following norm in $\M_+$:
$$
\|f\|_{\M_+}\df\inf\|\mu\|_{\M(\T)},
$$
the infimum being taken over all measures $\mu$ satisfying 
$\f(\z)=\int_\dd(1-\bar{\t}\z)^{-1}d\mu(\t)$, $\z\in\dd$. Similarly, we can define the space
$L^1_+$ of Cauchy integrals of $L^1$ functions,
$$
L^1_+\df\left\{\f:~\f(\z)=\int\limits_\dd\frac{h(\t)}{1-\bar{\t}\z}d\m(\t),~\z\in\dd,
~g\in L^1\right\}
$$
with an obvious definition of the norm in $L^1_+$. Clearly, if 
$\f(\z)=\int_\dd(1-\bar{\t}\z)^{-1}d\mu(\t)$, $\z\in\dd$, then $\G_\f=\G_\mu$.

\begin{thm}
\label{iCm}
Let $p>2/3<p$. Then there exists an analytic function $\psi$ in $\dd$ such that
$\G_\psi\in\Mp$ but $\psi\notin\M_+$.
\end{thm}

We deduce Theorem \ref{iCm} from imbedding theorems for Besov spaces (see Theorem \ref{besh1}
below).

We need some well-known facts about the {\it Dirichlet kernels}
$D_n$, $n\ge1$, defined by
$$
D_n(\z)\df\sum\limits_{k=-n}^n\z^k=\frac{\z^{n+1}-\z^{-n}}{\z-1},\quad\z\in\C,\quad\z\ne0.
$$ 
Obviously, $\|D_n\|_2=\sqrt{2n+1}$ and $\|D_n\|_\be=2n+1$, and so 
$\int\limits_{\T}|D_n|^p d\m\le(2n+1)^{p-1}$
for any $p\ge2$. Next, it is easy to see that
$$
|D_n(\z)|\le\min\left\{2n+1,\frac2{|1-\z|}\right\},\quad\z\in\T,
$$
which implies 
$$
\int\limits_{\T}|D_n|^p d\m\le C(p)(2n+1)^{p-1},\quad 1<p<\be.
$$ 
It is also clear that 
$$
|D_n(e^{\text{i}t})|\ge\re D_n(e^{\text{i}t})\ge n+\frac12,\quad 
t\in\left(-\frac\pi{3n},\frac\pi{3n}\right),
$$
whence,
$$
\int\limits_{\T}|D_n|^p d\m\ge\frac13n^{p-1},\quad n\ge1,\quad p>0.
$$
Let $Q_n\df\frac1{2n+1}D_n^2$.

\begin{lem}
\label{sum} 
There exists a constant $C$ such that
$$ 
\sum\limits_{n\ge0}\frac1{2^{n+1}+1}\left|Q_{2^n}\left(e^{\frac{\rm{i}n}{2^n}}\z\right)\right|
\le C
$$
for any $\z\in\T$.
\end{lem}

\Pf Set $U_n\df\{\z\in\T:10|\z-e^{-\frac{\text{i}n}{2^n}}|<\frac n{2^n}\}$ for $n\ge1$
and $U_0\df\T$.
It is easy to see that $\sum\limits_{n\ge0}\chi_{U_n}\le C$. Consequently,
$$
\sum\limits_{n\ge0}\frac1{2^{n+1}+1}\chi_{U_n}(\z)|Q_{2^n}(e^{\frac{\text{i}n}{2^n}}\z)|\le C,
quad\z\in\T.
$$
It remains to note that 
$\frac1{2^{n+1}+1}|Q_{2^n}(e^{\frac{\text{i}n}{2^n}}z)|<\frac{400}{n^2}$
if $z\not\in U_n$, and $\sum\limits_{n\ge1}\frac1{n^2}<+\infty$.
$\bl$

Let $N$ be a positive integer. We define the analytic polynomials $\Phi_n^{(N)}$ by
$$
\Phi_n^{(N)}(z)\df\frac{z^{4^n}}{(2^{n+1}+1)^N}(D_{2^n}(z))^{N+1},\quad n\ge N,\quad\z\in\C.
$$
Clearly, $|\Phi_n^{(N)}(z)|\le|Q_{2^n}(z)|$ for $z\in\T$.

\begin{thm}
\label{l/h}
Let $f(\z)\df\sum\limits_{n\ge N}a_n\Phi_n^N(e^{\frac{\rm{i}n}{2^n}}\z)$ for $\z\in\dd$.
Then the following statement are equivalent:

{\rm(i)} $\sum\limits_{n\ge N}|a_n|<+\infty$;

{\rm(ii)} $f\in H^1$;

{\rm(iii)} $f\in L^1_+$;

{\rm(iv)} $f\in \M_+$.
\end{thm}

\Pf The implications (ii)$\Rightarrow$(iii) and (iii)$\Rightarrow$(iv) are obvious.
To prove that (i) implies (ii), it suffices to observe that
$$
\left\|\Phi_n^{(N)}\right\|_1=(2^{n+1}+1)^{-N}\int\limits_{\T}|D_{2^n}|^{N+1}d\m\le1.
$$
It remains to prove that (iv) implies (i). Let $f\in\M_+$. Clearly, 
$$
\left|\sum_{k\ge0}\widehat f(k)\bar\a_k\right|\le C\|\sum_{k\ge0}\a_k z^k\|_\be
$$
for any polynomial $\sum_{k\ge0}\a_k\z^k$.
Put 
$$
a_n^*\df\left\{\begin{array}{ll}\frac{a_n}{|a_n|},&a_n\not=0,\\[.2cm]0,&a_n=0.
\end{array}\right.
$$
Let 
$$
g_m\df\sum_{n=N}^m \frac{a_n^*}{2^{n+1}+1}\Phi_n^N(e^{\frac{\text{i}n}{2^n}}z).
$$
By Lemma \ref{sum}, $\|g_m\|_\infty\le C$. Consequently,
$$
\sum\limits_{n=N}^m|a_n|(2^{n+1}+1)^{-1}\|\Phi_n^N\|_2^2\le C,\quad m\ge N.
$$
It remains to observe that 
\bey
(2^{n+1}+1)^{-1}\|\Phi_n^N\|_2^2&=&
(2^{n+1}+1)^{-2N-1}\int\limits_{\T}|D_{2^n}|^{2N+2}d\m\\[.2cm]
&\ge&\frac13(2^{n+1}+1)^{-2N-1}2^{(2N+1)n}\ge C(N).\quad\bl
\eey

We need the following well-known lemma.

\begin{lem} 
\label{p-1}
Let $f$ be a polynomial of degree at most $n$.
Then 
$$
\|f\|_\infty\le e(n+1)^{\1/p}\|f\|_p
$$
for any $p>0$.
\end{lem}

\Pf Let $n\ge1$. Put $g(z)\df f\left((1+\frac1n)z\right)$. Clearly, 
$$
\|g\|_p\le\left(1+\frac1n\right)^n\|f\|_p\le e\|f\|_p,\quad p>0.
$$
Using the well-known
inequality $|g(a)|\le\frac{\|g\|_p}{\left(1-|a|^2\right)^{1/p}}$, we obtain
$$
\|f\|_\infty\le\frac{e\|f\|_p}{\left(1-(1+\frac1n)^{-2}\right)^{1/p}}
\le e(n+1)^{1/p}\|f\|_p.\quad\bl
$$

\begin{cor} 
\label{p<1} 
Let $f$ be a polynomial of degree at most $n$.
Then 
$$
\|f\|_1\le e^{1-p}(n+1)^{1/{p_\#}}\|f\|_p
$$
for any $p\in(0,1)$.
$\bl$
\end{cor}

The following result is possibly known but we were unable to find a reference. 

\begin{thm}
\label{besh1}
Let $0<p\le1$ and $q>0$. Then the following statement are equivalent:

{\em(i)} $q\le1$;

{\em(ii)} $\pp_+(B_{p\,q}^{1/{p_\#}})\subset H^1$;

{\em(iii)} $\pp_+(B_{p\,q}^{1/{p_\#}})\subset L^1_+$;

{\em(iv)} $\pp_+(B_{p\,q}^{1/{p_\#}})\subset \M_+$.
\end{thm}

\Pf Let us prove that (i) implies (ii). Let $f\in\pp_+\left(B_{p\,q}^{1/{p_\#}}\right)$.
Then 
$$
\sum\limits_{n\ge0}2^{nq/p_\#}\|f*V_n\|_p^q<+\infty,
$$
where the $V_n$ are the polynomials defined in \S2. Consequently,
$$
\sum\limits_{n\ge0}2^{n/p_\#}\|f*V_n\|_p<+\infty,
$$
and by Corollary \ref{p<1},
$\sum\limits_{n\ge0}\|f*V_n\|_1<+\infty$. Thus, $f=\sum\limits_{n\ge0}f*V_n\in H^1$.

The implications (ii)$\Rightarrow$(iii) and (iii)$\Rightarrow$(iv) are trivial.
It remains to prove that (iv) implies (i). Take a positive integer $N$ such that
$p(N+1)>1$. Let $f$ be the function defined in the statement of Theorem \ref{l/h}. Clearly,
$f\in B_{p\,q}^{1/{p_\#}}$ if and only if 
$$
\sum\limits_{n\ge N} 2^{nq/p_\#}|a_n|^q\|\Phi_n^{(N)}\|_p^q<+\infty.
$$
Note that 
\bey
\int\limits_\T|\Phi_n^{(N)}|^pd\m&=&(2^{n+1}+1)^{-Np}\int\limits_{\T}|D_{2^n}|^{(N+1)p}d\m\\[.2cm]
&\le&C(p,N)2^{-Npn}2^{Npn+pn-n}=C(p,N)2^{-np/p_\#}.
\eey
Next,
$$
\int\limits_{\T}|\Phi_n^{(N)}|^pd\m=(2^{n+1}+1)^{-Np}\int\limits_{\T}|D_{2^n}|^{(N+1)p}d\m\ge
C(p,N)2^{-np/p_\#}.
$$ 
Consequently,
$f\in B_{p\,q}^{1/{p_\#}}$ if and only if 
$\sum\limits_{n\ge N}|a_n|^q<+\infty$. By Theorem \ref{l/h}, 
$f\in \M_+$ if and only if $\sum\limits_{n\ge N}|a_n|<+\infty$. Now it is obvious
that (iv) implies (i).
$\bl$

{\bf Proof of Theorem \ref{iCm}.} Let $p>\frac23$. 
Clearly, we can assume that $p\in(\frac23,1]$.
We have $p_\flat>1$, and by Theorem \ref{besh1}, there exists
a function $\psi\in \pp_+\left(B^{1/{p_\#}}_{p\,p_\flat}\right)\setminus \M_+$. It remains
to apply Theorem \ref{t3.y}. $\bl$

We do not know whether there exists a function
analytic in $\dd$ such that
$\Gamma_\psi\in\Mp$ for $p\le\frac23$ and
$\psi\not\in\M_+$.

\

\

\noindent
\begin{tabular}{p{8cm}p{14cm}}
A.B. Aleksandrov & V.V. Peller \\
St-Petersburg Branch & Department of Mathematics \\
Steklov Institute of Mathematics  & Kansas State University \\
Fontanka 27, 191011 St-Petersburg & Manhattan, Kansas 66506\\
Russia&USA
\end{tabular}

\end{document}